\documentclass[12pt]{article}
\usepackage{amsmath,amsthm, amssymb}
\def\Omega{\varOmega}
\def\ep{\varepsilon}

\def\Xint#1{\mathchoice
    {\XXint\displaystyle\textstyle{#1}}%
     {\XXint\textstyle\scriptstyle{#1}}%
     {\XXint\scriptstyle\scriptscriptstyle{#1}}%
     {\XXint\scriptstyle\scriptscriptstyle{#1}}%
	\!\int}
\def\XXint#1#2#3{{\setbox0=\hbox{$#1{#2#3}{\int}$}
	\vcenter{\hbox{$#2#3$}}\kern-.5\wd0}}

\begin{document}

\title{\large 
\bf The Harnack inequality and related properties \\
for solutions to elliptic and parabolic equations \\
with divergence-free lower-order coefficients}
\author { {\it A.I.~Nazarov, N.N.~Ural'tseva}\footnote {Partially supported
by RFBR grant 08-01-00748 and by grant NSh.4210.2010.1.},\\ 
Saint-Petersburg University, \\
{\small e-mail: \ al.il.nazarov@gmail.com}} 
\date{}\maketitle

\hfill To the memory of M.S. Birman

\section{Introduction}

Qualitative properties of solutions to partial differential equations are 
intensively studied over last half of century. In this paper we deal with
classical properties, namely, strong maximum principle, H\"older estimates, 
the Harnack inequality and the Liouville Theorem.\medskip

We consider elliptic and parabolic equations of divergence type:
$${\cal L}u \equiv -
D_i\big(a_{ij}(x)D_ju\big)+
b_i(x)D_iu=0;\eqno({\bf DE})$$
$${\cal M}u\equiv \partial_tu-
D_i\big(a_{ij}(x;t)D_ju\big)+
b_i(x;t)D_iu=0.\eqno({\bf DP})$$
We mostly deal with {\it a priori} estimates for Lipschitz generalized (sub/super)solutions.
When these estimates are established, we discuss the possibility of their generalization for weak (sub/super)solutions. In this case we assume 
$Du\in L_{2,loc}(\Omega)$ in ({\bf DE}) and $u\in L_{2,\infty,loc}(Q)$, 
$Du\in L_{2,loc}(Q)$ in ({\bf DP}). 

We always suppose that operators under consideration are uniformly elliptic 
(parabolic), i.e. for all values of arguments 
\begin{equation}\label{ell}
\nu|\xi|^2\le a_{ij}(\cdot)\xi_i\xi_j\le\nu^{-1}|\xi|^2,\qquad 
\xi\in{\mathbb R}^n,
\end{equation}
where $\nu$ is a positive constant.\medskip

The properties of generalized solutions to the equations ({\bf DE})--({\bf DP})
were investigated in a number of papers. H\"older estimates for solutions of 
({\bf DE}) were obtained by E.~De Giorgi \cite{DG} for ${\bf b}\equiv0$ and by 
C.~Morrey \cite{M} for ${\bf b}$ belonging to the Morrey space lying between $L_n$ 
and any $L_q$, $q>n$ (${\bf b}$ stands for $(b_i)$). Corresponding result for 
({\bf DP}) was established by J.~Nash \cite{Na} for ${\bf b}\equiv0$ and by 
O.A.~Ladyzhenskaya and N.N.~Ural'tseva \cite{LU1} for ${\bf b}\in L_{q+2}$, $q>n$.

Harnack's inequality for operators without lower-order coefficients was proved 
by J.~Moser (\cite{Mo1} for ({\bf DE}) and \cite{Mo2} for ({\bf DP})). 
N.~Trudinger \cite{Tru} proved it for ({\bf DE}) with 
${\bf b}\in L_q$, $q>n$. G.~Lieberman (see \cite[Ch. VI]{Li1}) extended the result of \cite{Mo2} 
for ${\bf b}\in L_{q,\ell}$, $\frac nq+\frac 2\ell<1$. Obviously, Harnack's inequality
implies H\"older estimates. Also strong maximum principle follows from 
Harnack's inequality and weak maximum principle. Some sharpening of mentioned
results, as well as corresponding results for nondivergence equations, are
discussed in our preprint \cite{NU}.

In this paper we consider mainly the equations ({\bf DE}) and ({\bf DP}) with additional
structure condition 
\begin{equation}\label{bezdiv}
\mbox{\rm div} ({\bf b})\le0\quad \mbox{in the sense of distributions}.
\end{equation}
The equations with the lower-order coefficients satisfying this structure condition 
arise in some applications (see, e.g., \cite{Z}, \cite{KNSS}, \cite{CSTY}, \cite{SSSZ}). 
We show that in this case the assumptions on 
${\bf b}$ can be considerably weakened in the scale of Morrey spaces.

Our paper is organized as follows. In Section 2 we deal with elliptic equations. 
Section 3 is devoted to parabolic equations (recall that only two-sided Liouville's Theorem holds
for these equations). In Section 4 we show an application of our results to some equations arising in
hydrodynamics. We underline that this Section contains just exemplary instances, and we make no pretence
to the novelty of results. In particular, the statements of Theorems 4.1 and 4.3 are in fact obtained in
\cite{KNSS}.\medskip

Let us recall some notation. $x=(x_1,\dots,x_n)$ is a vector 
in $\mathbb R^n$, $n\ge2$, with the Euclidean norm $|x|$; $(x;t)$ is a point 
in $\mathbb R^{n+1}$. 

$\Omega$ is a domain in $\mathbb R^n$ and 
$\partial\Omega$ is its boundary%; ${\text{\bf n}}=({\text{\bf n}}_i(x))$ is 
%the unit vector of the outward normal to $\partial\Omega$ at the point $x$
. For a cylinder $Q=\Omega\times]0,T[$ we denote by
$\partial''Q=\partial\Omega\times]0,T[$ its lateral surface and by
$\partial'Q=\partial''Q\cup\{\overline{\Omega} \times \{0\}\}$ its
parabolic boundary.

We define
$$\begin{array}{lll}
B_{R}(x^0)=\{x\in \mathbb R^n : |x-x^0|<R\}, & B_R=B_R(0);\\
Q_{R}^{\lambda,\theta}(x^0;t^0)=B_{\lambda R}(x^0)\times]t^0-\theta R^2;t^0[, & 
Q_{R}^{\lambda,\theta}=Q_{R}^{\lambda,\theta}(0;0), & Q_{R}=Q_{R}^{1,1}
\end{array}
$$
(note that $Q_{\lambda R}=Q_{R}^{\lambda,\lambda^2}$).

The indices $i, j$  vary from 1 to $n$.
Repeated indices indicate summation. 

The symbol $D_i$ denotes the operator of differentiation with respect to 
$x_i$; in particular, $Du=(D_1u, \dots, D_nu)$ is the gradient of 
$u$. $\partial_tu$ stands for the derivative of $u$ with respect to $t$.

The dashed integral stands for the mean value: 
$\Xint-\limits_Eu=(\mbox{meas}\, E)^{-1}\int\limits_Eu$.

We denote by $\|\cdot\|_{p,\Omega}$ the norm in $L_p(\Omega)$.
We introduce a scale of anisotropic spaces
$L_{q,\ell}(Q)=L_\ell\big(\,]0,T[\,\to L_q(\Omega)\big)$ with the norm
$\|f\|_{q,\ell,Q}=\big\|\|f(\cdot;t)\|_{q,\Omega}\big\|_{\ell,]0,T[}$.
\noindent Obviously, $L_{q,q}(Q)=L_q(Q)$.

We also introduce a scale of Morrey spaces
$$\mathbb M^{\alpha}_q(\Omega)=\{f\in L_q(\Omega):\ \|f\|_{\mathbb M^{\alpha}_q(\Omega)}\equiv
\sup\limits_{B_\rho(x)\subset\Omega}\rho^{-\alpha}\|f\|_{q,B_\rho(x)}<\infty\}.$$
Parabolic Morrey spaces $\mathbb M^{\alpha}_{q,\ell}(Q)$ 
are introduced in a similar way, using $Q_\rho(x;t)\subset Q$ instead of 
$B_\rho(x)\subset\Omega$.

Finally, we introduce the space ${\cal V}(Q)$ of weak solutions to ({\bf DP}) with the norm defined by
$$\|f\|^2_{{\cal V}(Q)}=\|f\|^2_{2,\infty,Q}+\|Df\|^2_{2,2,Q}.
$$

We set $f_+=\max \{f,0\},\ \ f_-=\max \{-f,0\}$, \
${\text{osc}}_{\Omega}  f =\sup_{\Omega}  f - \inf_{\Omega}  f$. For $1\le p<n$,
$p^*=\frac {np}{n-p}$ is the Sobolev conjugate to $p$.

We use letters $N$, $C$ (with or without indices) to denote various 
constants. To indicate that, say, $N$ depends on some parameters, we list 
them in the parentheses: $N(\dots)$.

\section{Elliptic case}

Recall that $u$ is a (Lipschitz) subsolution of the equation ${\cal L}u=0$ in $\Omega$ (here $\cal L$ 
is an operator of the form (\,{\bf DE})), if for any Lipschitz test function $\eta\ge0$, supported in $\Omega$, 
$$\int\limits_{\Omega}(a_{ij}D_juD_i\eta+b_iD_iu\,\eta)\,dx\le0.
$$
We take $\eta=\varphi'(u)\cdot\xi$, where $\xi$ is a nonnegative Lipschitz function, supported in 
$B_{\lambda R}\subset\Omega$, while $\varphi\in{\cal C}^2(\mathbb R)$ is a convex function vanishing in 
$\mathbb R_-$. This gives
\begin{equation}\label{Moser}
\int\limits_{B_{\lambda R}\cap\{u>0\}}\Big(a_{ij}D_jvD_i\xi+\frac {\varphi''(u)}
{\varphi^{\prime2}(u)}\,a_{ij}D_jvD_iv\,\xi+b_iD_iv\,\xi\Big)\,dx\le0,
\end{equation}
ЦДЕ $v=\varphi(u)$.

Then, by mollification at a neighborhood of the origin, one can weaken in (\ref{Moser}) the assumption
$\varphi\in{\cal C}^2(\mathbb R)$ to $\varphi\in{\cal C}^2(\mathbb R_+\cup\mathbb R_-)$.\medskip

{\bf Lemma 2.1}. {\it Let $\cal L$ be an operator 
of the form (\,{\bf DE}) in $B_{\lambda R}(x^0)$, $\lambda>1$,  and let the conditions (\ref{ell}) and 
(\ref{bezdiv}) be satisfied. Let also
${\bf b}\in L_q(B_{\lambda R}(x^0))$
with some $\frac n2<q\le n$\footnote{For $q=n$, the assumption (\ref{bezdiv}) can be removed. We discuss this at the end of this Section.}.

Then there exists a positive constant $N_1$ depending on $n$, $\nu$, $\lambda$, $q$ and 
the quantity 
$${\cal N}={\cal N}(R,\lambda)\equiv R^{1-\frac nq}\|{\bf b}\|_{q,B_{\lambda R}(x^0)},$$
such that any Lipschitz subsolution 
of the equation ${\cal L}u=0$ in $B_{\lambda R}(x^0)$ satisfies}
\begin{equation}\label{estmax}
\sup\limits_{B_R(x^0)} u_+\le N_1
\bigg(\!\!\Xint{\quad\ \,-}\limits_{B_{\lambda R}(x^0)}\!\!u_+^2dx\bigg)^{\frac 12}.
\end{equation}

\begin{proof} We use classical technique of Moser (see, e.g., \cite[Ch.IX]{LU2}). 
Without loss of generality, we assume $x^0=0$. 

We put in (\ref{Moser}) $\varphi(\tau)=\tau_+^p$, $p>1$, and $\xi=v\zeta^2$ 
where $\zeta$ is a smooth cut-off function in $B_{\lambda R}$. Then we obtain
\begin{equation}\label{Moser_power}
\int\limits_{B_{\lambda R}}\Big(\frac {2p-1}{p}\,a_{ij}D_jvD_iv\,\zeta^2
+2a_{ij}D_jv\,vD_i\zeta\,\zeta+b_iD_iv\,v\zeta^2\Big)\,dx\le0.
\end{equation}
The last term in (\ref{Moser_power}) can be estimated using (\ref{bezdiv}) and the H\"older
inequality:
\begin{equation}\label{bezdiv1}
-\int\limits_{B_{\lambda R}}b_iD_iv\,v\zeta^2\,dx\le
\int\limits_{B_{\lambda R}}b_iv^2\zeta\,D_i\zeta\,dx
\le\|{\bf b}\|_{q,B_{\lambda R}}\|v\zeta\|^{2-\frac 1s}_{r,B_{\lambda R}}
\|v\zeta^{1-s}|D\zeta|^s\|^{\frac 1s}_{2,B_{\lambda R}},
\end{equation}
where $s>2$ is defined by $\frac 1s=1-\frac n{2q}$ while $r=\frac {2(2q+n)}{2q+n-4}$. Note that $2<r<2^*$, and, by the embedding theorem,
\begin{equation}\label{embed}
\|v\zeta\|_{r,B_{\lambda R}}\le
C(n)(\lambda R)^{n(\frac 1r-\frac 1{2^*})}\,\|D(v\zeta)\|_{2,B_{\lambda R}}
\le C(n)(\lambda R)^{\frac 1{2s-1}}\,\Big(\|Dv\,\zeta\|_{2,B_{\lambda R}}+\|vD\zeta\|_{2,B_{\lambda R}}\Big).
\end{equation}

Using (\ref{ell}), (\ref{bezdiv1}) and (\ref{embed}), we obtain from (\ref{Moser_power})
\begin{multline*}%\label{aa}
\|Dv\,\zeta\|^2_{2,B_{\lambda R}}\le
\frac 1\nu\int\limits_{B_{\lambda R}}a_{ij}D_jvD_iv\,\zeta^2\,dx\le C_1(n,\nu,s,\lambda)\times\\
\times\bigg[\|Dv\,\zeta\|_{2,B_{\lambda R}}\|vD\zeta\|_{2,B_{\lambda R}}+
R^{\frac 1s}\|{\bf b}\|_{q,B_{\lambda R}}\Big(\|Dv\,\zeta\|^{2-\frac 1s}_{2,B_{\lambda R}}+
\|vD\zeta\|^{2-\frac 1s}_{2,B_{\lambda R}}\Big)
\|v\zeta^{1-s}|D\zeta|^s\|^{\frac 1s}_{2,B_{\lambda R}}\bigg],
\end{multline*}
and therefore
\begin{equation}\label{ee}
\|Dv\,\zeta\|_{2,B_{\lambda R}}\le
C_2(n,\nu,s,\lambda)\cdot\bigg[\|vD\zeta\|_{2,B_{\lambda R}}+R\|{\bf b}\|^s_{q,B_{\lambda R}}
\|v\zeta^{1-s}|D\zeta|^s\|_{2,B_{\lambda R}}\bigg].
\end{equation}

We put $R_m=R(1+2^{-m}(\lambda-1))$, $m\in \mathbb N\cup\{0\}$, and substitute $\zeta=\zeta_m$ such that 
$$ \zeta_m\equiv1\ \ \mbox{in}\ \ B_{R_{m+1}};\quad \zeta_m\equiv0\ \ \mbox{out of}\ \ B_{R_m};
\qquad \frac{|D\zeta_m|}{\zeta_m^{1-\frac 1s}}\le \frac {2^mC_3(s)}{(\lambda-1)R}.$$ 
Then (\ref{ee}) implies
\begin{equation}\label{eee}
\|Dv\,\zeta_m\|_{2,B_{R_m}}\le
\frac {C_4(n,\nu,s,\lambda)}R\cdot\|v\|_{2,B_{R_m}}\cdot
\big(2^m+\big(2^m{\cal N}\big)^s\big).
\end{equation}

Now for $p=p_m\equiv (\frac r2)^m$
 we obtain from (\ref{embed}) and (\ref{eee})
\begin{multline}\label{iteration}
\bigg(\Xint{\quad\ \, -}\limits_{B_{R_{m+1}}}\!\! u_+^{2p_{m+1}}dx\bigg)^{\frac {1}{2p_{m+1}}}\!\!\le
\bigg(C(n)\Xint{\ \ -}\limits_{B_{R_m}} (v\zeta_m)^r dx\bigg)^{\frac {1}{rp_m}}\le\\
\le 
\bigg(2^{2ms}C_5\Xint{\ \ -}\limits_{B_{R_m}} v^2\,dx\bigg)
^{\frac {1}{2p_m}}\!\!=
\bigg(2^{2ms}C_5\Xint{\ \ -}\limits_{B_{R_m}} u_+^{2p_m}
dx\bigg)^{\frac 1{2p_m}}\!,
\end{multline}
where $C_5$ depends only on $n$, $\nu$, $\lambda$, $s$ and ${\cal N}$.

Iterating (\ref{iteration}) we arrive at (\ref{estmax}).
\end{proof}

{\bf Corollary 2.1}. {\it Let $\cal L$ satisfy the assumptions of Lemma 2.1 in $B_{\lambda R}(x^0)$. 
If a Lipschitz subsolution of ${\cal L}u=0$ in $B_{\lambda R}(x^0)$ satisfies}
\begin{equation}\label{tiny}
\mbox{meas}\, (\{u> k\}\cap B_{\lambda R}(x^0))\le \mu\,\mbox{meas}\,(B_{\lambda R}),
\qquad \mu< N_1^{-2},
\end{equation}
{\it for some $k$, then}
\begin{equation}\label{estmax1}
\sup\limits_{B_R(x^0)} (u-k)\le N_1\sqrt{\mu}\sup\limits_{B_{\lambda R}(x^0)}(u-k),
\end{equation}
{\it (here $N_1$ is the constant from Lemma 2.1).}

\begin{proof} We apply Lemma 2.1 to $u-k$.
\end{proof}

We need the following variant of the embedding theorem.\medskip

{\bf Proposition A}. {\it Let $1\le p<n$. Suppose that a non-negative function $u\in W^1_p(B_R)$ vanishes 
on a positive measure set ${\cal E}_0$. Let $\eta=\eta(|x|)$ be a non-decreasing function, $0\le\eta\le 1$,
and $\eta\big|_{{\cal E}_0}\equiv 1$. Then, for any $1\le q\le p^*$ and for any measurable set
${\cal E}\subset B_R$,}
$$\|u\,\eta\|_{q,{\cal E}}\le \frac {C(n)R^n}{\mbox{meas}\,({\cal E}_0)}\ 
\mbox{meas}^{\frac 1q-\frac 1{p^*}}({\cal E})\cdot\|Du\,\eta\|_{p,B_R}.
$$

\begin{proof}
For $q=p=1$ this Proposition was proved in \cite[Ch. II, Lemma 5.1]{LSU}. In this Lemma the following
inequality was obtained:
$$\mbox{meas}\,({\cal E}_0)\cdot u(x)\,\eta(x)\le
\frac {(2R)^n}{n}\int\limits_{B_R}\frac {|Du(y)|\,\eta(y)}{|y-x|^{n-1}}\,dy.
$$
By the Hardy--Littlewood--Sobolev inequality (see, e.g., \cite[Sec. 4.3]{LL}), we get
$$\mbox{meas}\,({\cal E}_0)\cdot \|u\,\eta\|_{p^*,B_R}\le C(n,p)R^n\cdot \|Du\,\eta\|_{p,B_R},
$$
and the statement follows by H\"older inequality.
\end{proof}

{\bf Lemma 2.2}. {\it Let $\cal L$ satisfy the assumptions of Lemma 2.1 in $B_{\lambda R}(x^0)$. 
Then for any $\delta>0$ there exists a positive constant $\beta$ depending on $n$, $\nu$, $\lambda$, 
$q$, $\delta$ and the quantity ${\cal N}$, such that if 
a Lipschitz nonnegative supersolution of ${\cal L}V=0$ in $B_{\lambda R}(x^0)$ satisfies}
\begin{equation}\label{tiny1}
\mbox{meas}\, (\{V\ge k\}\cap B_R(x^0))\ge\delta\cdot \mbox{meas}\, (B_R) 
\end{equation}
{\it for some $k>0$, then}
\begin{equation}\label{estmin}
\inf\limits_{B_R(x^0)} V\ge \beta k.
\end{equation}

\begin{proof} Without loss of generality, we can assume $V>0$; otherwise we deal with $V+\ep$ and pass
to the limit as $\ep\downarrow0$. Also we put $x^0=0$. 

We define $u=1-\frac Vk$. Note that $u<1$ is a subsolution, and therefore, we can apply the relation 
(\ref{Moser}) with $\varphi$ defined only for $\tau<1$.

We put in (\ref{Moser}) $\varphi(\tau)=\ln_-(1-\tau)$. This gives for $v=\varphi(u)$
\begin{equation}\label{Moser_ln}
\int\limits_{B_{\lambda R}}\Big(a_{ij}D_jvD_i\xi+a_{ij}D_jvD_iv\,\xi+b_iD_iv\,\xi\Big)\,dx\le0.
\end{equation}
We substitute into (\ref{Moser_ln}) $\xi=\zeta^2$ where $\zeta$ is a smooth cut-off function that 
equals $1$ in $B_{\frac{1+\lambda}2R}$. Then, using (\ref{ell}), (\ref{bezdiv}) and the H\"older inequality, 
we obtain
\begin{multline*}%\label{aa}
\|Dv\,\zeta\|^2_{2,B_{\lambda R}}\le
\frac 1\nu\int\limits_{B_{\lambda R}}a_{ij}D_jvD_iv\,\zeta^2\,dx\le \frac 2\nu
\int\limits_{B_{\lambda R}}\Big(-a_{ij}D_jv\zeta D_i\zeta+b_iv\zeta\,D_i\zeta\Big)\,dx\le\\
\le \frac 2\nu (\|Dv\,\zeta\|_{2,B_{\lambda R}}\|D\zeta\|_{2,B_{\lambda R}}+
\|{\bf b}\|_{q,B_{\lambda R}}\|v\,\zeta\|_{q',B_{\lambda R}}\|D\zeta\|_{\infty,B_{\lambda R}}).
\end{multline*}
Note that $v$ vanishes on the set $\{V\ge k\}\cap B_R$. Therefore, we can estimate the last term by 
Proposition A. By (\ref{tiny1}), this gives
$$\|Dv\,\zeta\|_{2,B_{\lambda R}}\le C_6(n,\nu,\lambda,q,\delta)R^{\frac n2-1}\cdot
\big(1+{\cal N}\big).
$$
Applying Proposition A once more, we obtain
$$\bigg(\!\!\Xint{\quad\ \,-}\limits_{B_{\frac {1+\lambda}2 R}}\!\!v^2dx\bigg)^{\frac 12}\le
C_7,$$
where $C_7$ depends only on $n$, $\nu$, $\lambda$, $q$, $\delta$ and ${\cal N}$.

Finally, the relation (\ref{Moser_ln}) implies that $v$ is a subsolution. So, we apply Lemma 2.1 to $v$ in
$B_{\frac {1+\lambda}2 R}$ and arrive at the estimate $\sup\limits_{B_R} v_+\le C_8\equiv N_1C_7$, 
which is equivalent to (\ref{estmin}) with $\beta=\exp(-C_8)$.
\end{proof}

{\bf Corollary 2.2} (strong maximum principle). {\it Let $\cal L$ be an operator of the form (\,{\bf DE}) 
in $\Omega$,  and let the conditions (\ref{ell}) and (\ref{bezdiv}) be satisfied. Let also
${\bf b}\in L_{q,loc}(\Omega)$ with some $\frac n2<q\le n$. Then any Lipschitz nonconstant supersolution of 
${\cal L}V=0$ in $\Omega$ cannot attain its minimum at interior point of $\Omega$.}\medskip

\begin{proof} Assume the converse. Without loss of generality, $\inf\limits_\Omega V=0$.
Then there exists $x^0\in\Omega$ which is a frontier 
point of the set $\{V>0\}$. Choose $R$ such that $\overline{B_{2R}}(x^0)\subset\Omega$. 
Then the relation (\ref{tiny1}) holds for some $k>0$ and $\delta>0$, and we obtain (\ref{estmin}),
a contradiction.
\end{proof}

{\bf Lemma 2.3}. {\it Let $\cal L$ satisfy the assumptions of Lemma 2.1 in $B_{3R}$. 
Then any Lipschitz solution of ${\cal L}u=0$ in $B_{3R}$ satisfies the estimate}
\begin{equation}\label{osc}
\underset{B_R}{\mbox{osc}}\ u\le \varkappa_0\,\underset{B_{3R}}{\mbox{osc}}\ u,
\end{equation}
{\it where $\varkappa_0<1$ depends on $n$, $\nu$, $q$ and the quantity ${\cal N}$.}

\begin{proof} We set 
$$k=\frac 12\big(\sup\limits_{B_{3R}}u+\inf\limits_{B_{3R}}u\big)$$
and consider two cases.\medskip

{\bf 1}. Let the relation (\ref{tiny}) hold with $\lambda=2$ and $\mu=\frac 14N_1^{-2}$. 
Then, by Corollary 2.1,
$$\sup\limits_{B_R} u\le \frac 12 \big(\sup\limits_{B_{2R}}u+k\big)\le
\sup\limits_{B_{3R}}u-\frac 14\,\underset{B_{3R}}{\mbox{osc}}\ u.
$$

{\bf 2}. In the opposite case we apply Lemma 2.2 (with $\lambda=\frac 32$ and $\delta=\mu$) to the (non-negative) function $V=u-\inf\limits_{B_{3R}}u$. This gives
$\inf\limits_{B_R} V\ge\inf\limits_{B_{2R}} V\ge \beta\big(k-\inf\limits_{B_{3R}}u\big)$, and thus,
$$\inf\limits_{B_R} u\ge\inf\limits_{B_{3R}} u+\frac {\beta}2\,
\underset{B_{3R}}{\mbox{osc}}\ u.
$$

In both cases we arrive at (\ref{osc}) with $\varkappa_0=\min\big\{\frac 14,\frac{\beta}2\big\}$.
\end{proof}

{\bf Corollary 2.3} (H\"older estimate). {\it Let $\cal L$ satisfy the assumptions of Lemma 2.1
in $B_{R_0}$. Let also $\sup\limits_{R<R_0}{\cal N}(R,1)<\infty$.

Then any Lipschitz solution of ${\cal L}u=0$ in $B_{R_0}$ satisfies the estimate}
\begin{equation}\label{Holder}
 \underset{B_\rho}{\mbox{osc}}\ u\le N_2 \Big(\frac \rho r\Big)^\gamma
\cdot\underset{B_r}{\mbox{osc}}\ u,\qquad 0<\rho<r\le R_0,
\end{equation}
{\it where $N_2$ and $\gamma$ depend on $n$, $\nu$, $q$ and $\sup\limits_{R<R_0}{\cal N}(R,1)$.}

\begin{proof} 
Iterating the estimate (\ref{osc}) we arrive at (\ref{Holder}) with $\gamma=-\log_3(\varkappa_0)$.
\end{proof}

{\bf Corollary 2.4} (two-sided Liouville's theorem). {\it Let $\cal L$ be an operator 
of the form (\,{\bf DE}) in $\mathbb R^n$, and let the conditions (\ref{ell}) and 
(\ref{bezdiv}) be satisfied. Let also ${\bf b}\in L_{q,loc}(\mathbb R^n)$, with some 
$\frac n2<q\le n$. Finally, assume that}
\begin{equation}\label{Liouville}
\liminf\limits_{R\to\infty}{\widehat{\cal N}}(R,1)<\infty.
\end{equation}
{\it Then any Lipschitz bounded solution of ${\cal L}u=0$ in $\mathbb R^n$ is a 
constant.}\medskip

{\bf Remark 1}. If ${\bf b}\in L_q(\mathbb R^n)$, then (\ref{Liouville}) is 
obviously satisfied.

\begin{proof} Iteration of (\ref{osc}) with respect to a suitable sequence $R_m\to\infty$ gives 
the statement.
\end{proof}

{\bf Lemma 2.4}. {\it Let $\cal L$ be an operator of the form (\,{\bf DE}) 
in $B_{2R}$,  and let the conditions (\ref{ell}) and (\ref{bezdiv}) be satisfied. Let also
${\bf b}\in \mathbb M^{\frac nq-1}_q(B_{2R})$ with some $\frac n2<q\le n$.

Let for a Lipschitz nonnegative supersolution of ${\cal L}V=0$ in $B_{2R}$ and for some $y\in B_{2R}$, the 
inequality $\inf\limits_{B_\rho(y)}V=k>0$ holds with $\rho=\frac 14\mbox{\rm dist}(y,\partial B_{2R})$.
Then}
\begin{equation}\label{estmin1}
\inf\limits_{B_R} V\ge \widehat\beta \Big(\frac \rho R\Big)^{\widehat\gamma} k,
\end{equation}
{\it where $\widehat\beta $ and $\widehat\gamma$ depend on $n$, $\nu$, $q$ and 
$\|{\bf b}\|_{\mathbb M^{\frac nq-1}_q(B_{2R})}$.}

\begin{proof} Denote by ${\mathfrak N}$ an integer number such that 
$2^{-({\mathfrak N}+1)}R<\rho\le 2^{-{\mathfrak N}}R$ and consider a ball $B_{{\mathfrak r}_0}(y^0)$, where 
${\mathfrak r}_0=2^{-{\mathfrak N}}R$, $y^0=2R(1-2^{-{\mathfrak N}}){\bf e}$ and ${\bf e}=\frac y{|y|}$. It is
easy to see that $B_{{\mathfrak r}_0}(y^0)\subset B_{3\rho}(y)$, and by Lemma 2.2 (with $\lambda=\frac 43$ 
and $\delta=\frac 1{3^n}$), 
\begin{equation}\label{ff}
\inf\limits_{B_{{\mathfrak r}_0}(y^0)} V\ge\inf\limits_{B_{3\rho}(y)} V\ge \beta k.
\end{equation}

Now we introduce the sequence of balls $B_{{\mathfrak r}_m}(y^m)$, $m=1,\dots,{\mathfrak N}$, as follows:
$${\mathfrak r}_m=2{\mathfrak r}_{m-1},\qquad y^m=y^{m-1}-{\mathfrak r}_m{\bf e}.$$
For all $m=1,\dots,{\mathfrak N}$ one has
$$B_{2{\mathfrak r}_m}(y^m)\subset B_{2R};\qquad\mbox{meas}\,(B_{{\mathfrak r}_{m-1}}(y^{m-1})\cap 
B_{{\mathfrak r}_m}(y^m))\ge C(n)\cdot\mbox{meas}\,(B_{{\mathfrak r}_m}).
$$
Thus, Lemma 2.2 (with $\lambda=2$ and $\delta=C(n)$) gives
$$\inf\limits_{B_{{\mathfrak r}_m}(y^m)} V\ge \beta\cdot\inf\limits_{B_{{\mathfrak r}_{m-1}}(y^{m-1})} V.$$

Since $B_{{\mathfrak r}_{\mathfrak N}}(y^{\mathfrak N})=B_R$, we obtain
$$\inf\limits_{B_R} V\ge \beta^{\mathfrak N}\cdot\inf\limits_{B_{{\mathfrak r}_0}(y^0)} V\ge 
\Big(\frac \rho R\Big)^{\widehat\gamma}\cdot\inf\limits_{B_{{\mathfrak r}_0}(y^0)} V,$$
where $\widehat\gamma=-\log_2(\beta)$.

Combining this estimate with (\ref{ff}), we arrive at (\ref{estmin1}).
\end{proof}

{\bf Theorem 2.5} (the Harnack inequality). {\it Let $\cal L$ satisfy the assumptions of Lemma 2.4 
in $B_{2R}$. Then there exists a positive constant $N_3$ depending on $n$, $\nu$, $q$
and $\|{\bf b}\|_{\mathbb M^{\frac nq-1}_q(B_{2R})}$, such that any Lipschitz nonnegative solution 
of ${\cal L}u=0$ in $B_{2R}$ satisfies}
\begin{equation}\label{Harnack}
\sup\limits_{B_R}u\le N_3\cdot\inf\limits_{B_R}u.
\end{equation}

\begin{proof} We follow the idea of Safonov (\cite{S2}). 
Denote by $y\in B_{2R}$ a maximum point of the function
$$v(x)=(\mbox{dist}(x,\partial B_{2R}))^{\widehat\gamma}\!\cdot u(x)$$
(here $\widehat\gamma$ is the constant from Lemma 2.4) and set
$$\rho=\frac 14\,\mbox{dist}(y,\partial B_{2R});\qquad {\mathfrak M}=v(y)=(4\rho)^{\widehat\gamma}\!\cdot u(y).
$$

It is obvious that
\begin{eqnarray}
\sup\limits_{B_R}u\le \frac {\mathfrak M}{R^{\widehat\gamma}}=
\Big(\frac {4\rho} R\Big)^{\widehat\gamma}\!\cdot u(y);\label{sup}\\
\sup\limits_{B_{2\rho}(y)}u\le \frac {\mathfrak M}{(2\rho)^{\widehat\gamma}}=
2^{\widehat\gamma}\cdot u(y).\phantom{\rho}\label{sup1}
\end{eqnarray}

Denote $k=\frac 12 u(y)$. If
$\mbox{meas}\, (\{u> k\}\cap B_{2\rho}(y))\le \mu\,\mbox{meas}\,(B_{2\rho})$,
then Corollary 2.1 (with $\lambda=2$) and (\ref{sup1}) imply the relation
$$k=u(y)-k\le \sup\limits_{B_\rho(y)}(u-k)\le N_1\sqrt{\mu}\sup\limits_{B_{2\rho}(y)}(u-k)
\le N_1\sqrt{\mu}\,(2^{\widehat\gamma+1}-1)k,
$$
which is impossible for $\mu\le\mu_0\equiv\frac 1{2^{2\widehat\gamma+2}}\,N_1^{-2}$. Thus,
$\mbox{meas}\, (\{u> k\}\cap B_{2\rho}(y))\ge \mu_0\,\mbox{meas}\,(B_{2\rho})$, and Lemma 2.2
(with $\lambda=2$ and $\delta=\mu_0$) gives
\begin{equation}\label{inf}
\inf\limits_{B_\rho(y)}u\ge\inf\limits_{B_{2\rho}(y)}u\ge \beta k=\frac{\beta}2\cdot u(y).
\end{equation}

Finally, Lemma 2.4 gives
\begin{equation}\label{inf1}
\inf\limits_{B_R}u\ge \widehat\beta \Big(\frac \rho R\Big)^{\widehat\gamma}\inf\limits_{B_\rho(y)}u.
\end{equation}
Combining (\ref{sup}), (\ref{inf}) and (\ref{inf1}), we arrive at (\ref{Harnack}) with 
$N_3=\frac{2^{2\widehat\gamma+1}}{\beta\widehat\beta }$.
\end{proof}

{\bf Theorem 2.6} (one-sided Liouville's theorem). {\it Let $\cal L$ be an operator 
of the form (\,{\bf DE}) in $\mathbb R^n$, and let the conditions (\ref{ell}) and 
(\ref{bezdiv}) be satisfied. Let also
${\bf b}\in \mathbb M^{\frac nq-1}_{q,loc}(\mathbb R^n)$ 
with some $\frac n2<q\le n$, and for some $\delta>0$}
\begin{equation}\label{Liouville1}
\liminf\limits_{R\to\infty}\sup\limits_{|x|=R}
\|{\bf b}\|_{\mathbb M^{\frac nq-1}_q(B_{\delta R}(x))}<\infty.
\end{equation}
{\it Then any Lipschitz semibounded solution of ${\cal L}u=0$ in $\mathbb R^n$ is a 
constant.}\medskip

{\bf Remark 2}. If ${\bf b}\in \mathbb M^{\frac nq-1}_q(\mathbb R^n)$, 
then (\ref{Liouville1}) is obviously satisfied.

\begin{proof} Without loss of generality, we can assume that $u$ is bounded from below, and 
$\inf\limits_{\mathbb R^n} u=0$.

We take a sequence $R_m\to\infty$ such that 
$\mathfrak B\equiv\sup\limits_m\sup\limits_{|x|=R_m}\|{\bf b}\|_{\mathbb M^{\frac nq-1}_q(B_{\delta R}(x))}<\infty$. 
Further, we cover the sphere $|x|=1$ with a finite set of balls 
$B_{\frac{\delta} 2}(x)$ and dilate these balls to the covering of the sphere $|x|=R_m$. 
Applying Theorem 2.5 to all the balls of this covering, we obtain 
$\sup\limits_{|x|=R_m}u\le C(n,\nu,q,\mathfrak B,\delta)\cdot\!
\inf\limits_{|x|=R_m}\!u$ for any $m$. By Corollary 2.2,
$$\sup\limits_{B_R}u=\sup\limits_{|x|=R}u;\qquad
\inf\limits_{|x|=R}u=\inf\limits_{B_{R}}u\to0\quad \mbox{as}\quad R\to\infty,$$ 
and the statement follows.
\end{proof}

Let us discuss briefly the possibility to generalize all previous statements for weak (sub/super)solutions.\medskip

The proof of Lemma 2.1 runs without changes\footnote{More formally, in this case the inequality (\ref{Moser})
holds under additional condition that $\varphi$ is globally Lipschitz. Thus, to derive (\ref{ee}) 
one should take $\varphi'(u_+)=p\,\min\{u_+,N\}^{p-1}$, $\xi=\min\{u_+,N\}^p\zeta^2$, $N>0$, 
and then pass to the limit as $N\to\infty$.} also for weak subsolutions of ${\cal L}u=0$ if the bilinear form 
$${\cal B}\big\langle u,\eta\big\rangle\equiv \int\limits_{B_{\lambda R}}b_iD_iu\,\eta\,dx$$
can be continuously extended to the pair $(v,v\zeta^2)$ with $Dv\in L_2(B_{\lambda R})$. This is certainly true provided
$$\big|{\cal B}\big\langle u,\eta\big\rangle\big|\le 
C\|Du\|_{2,B_{2\lambda R}}\|D\eta\|_{2,B_{2\lambda R}},\qquad u,\eta\in {\cal C}^\infty_0(B_{2\lambda R}).
$$
It is shown in \cite{MV} that the last estimate holds if 
\begin{equation}\label{mazya}
\Delta^{-1} {\rm rot}({\bf b})\in BMO^{n\times n}(B_{2\lambda R});\qquad
h=|\nabla(\Delta^{-1} {\rm div}({\bf b}))|^2\in{\mathfrak M}^{1,2}_+
\end{equation}
(here we assume ${\bf b}$ extended by zero); the last notation means a class of so-called admissible weights, i.e.
$$\int\limits_{B_{2\lambda R}}h|v|^2\,dx\le C\|Dv\|^2_{2,B_{2\lambda R}},\qquad 
v\in{\cal C}^\infty_0(B_{2\lambda R}).
$$

If ${\bf b}\in\mathbb M^{\frac nq-1}_q(B_{\lambda R})$, the first relation in (\ref{mazya}) follows from elliptic coercive estimates and the Poincar\'e inequality, see, e.g., \cite{Tro}. Thus, Lemma 2.1 and, therefore, all subsequent statements hold true for weak (sub/super)solutions of ${\cal L}u=0$ if, for example,
 ${\bf b}\in\mathbb M^{\frac nq-1}_q$ and ${\rm div}({\bf b})\equiv0$.\medskip

In addition, let us consider the case $q=n$. The space $\mathbb M^{\frac nq-1}_q(\Omega)$ now becomes 
conventional Lebesgue space $L_n(\Omega)$, and we claim that main results of this section hold true without 
the assumption (\ref{bezdiv})\footnote{Moreover, in this case the assumption ${\bf b}\in L_n$ can be weakened in the scale of Lorentz spaces to ${\bf b}\in\Lambda_{n,q}$ with any $q<\infty$. We do not discuss it here for the reason of place.}. Note that in this case the relation (\ref{Moser_power}) is fulfilled for any weak subsolution $u$. 

First, let $n\ge3$. Then we estimate the last term in (\ref{Moser_power}) by the H\"older inequality and the Sobolev
inequality and obtain an analog of (\ref{ee}):
$$\|Dv\,\zeta\|_{2,B_{\lambda R}}\le
C'_2(n,\nu)\cdot\Big[\|vD\zeta\|_{2,B_{\lambda R}}+
\|{\bf b}\|_{n,B_{\lambda R}}\|Dv\,\zeta\|_{2,B_{\lambda R}}\Big].
$$
If $\|{\bf b}\|_{n,B_{\lambda R}}\le\ep(n,\nu)\equiv (2C'_2)^{-1}$, then for $p=p_m\equiv 2\big(\frac n{n-2}\big)^m$
we obtain an analog of (\ref{iteration}):
$$\bigg(\Xint{\quad\ \, -}\limits_{B_{R_{m+1}}}\!\! u_+^{2p_{m+1}}dx\bigg)^{\frac {1}{2p_{m+1}}}\!\!\le
\bigg(C'_5\Xint{\ \ -}\limits_{B_{R_m}} u_+^{2p_m}dx\bigg)^{\frac 1{2p_m}}\!,
$$
where $C'_5$ depends only on $n$, $\nu$, and $\lambda$. The remainder of the proof of Lemma 2.1 runs without changes.

Similarly we prove Lemmas 2.2 and 2.4 for weak supersolutions of ${\cal L}V=0$ under the same assumption 
$\|{\bf b}\|_{n,B_{\lambda R}}\le\ep(n,\nu)$ (in Lemma 2.4 $\lambda=2$).

In the case $n=2$ we use the Yudovich--Pohozhaev embedding theorem (see, e.g., \cite[10.6]{BIN}) instead of 
the Sobolev inequality. This gives us Lemmas 2.1, 2.2, 2.4 under the assumption 
$\|{\bf b}\ln^{\frac 12}\big(1+\lambda R|{\bf b}|\big)\|_{2,B_{\lambda R}(x^0)}\le \ep(n,\nu)$.\medskip

Further, strong maximum principle holds without smallness assumptions on ${\bf b}$. Indeed, one 
can choose $R$ sufficiently small such that these assumptions are fulfilled.\medskip

Since the proof of Theorem 2.5 depends only on Lemmas 2.1, 2.2, 2.4, the Harnack inequality evidently holds under smallness assumption on ${\bf b}$. However, we can exclude this assumption using a trick of M.V. Safonov.\cite{S3}\medskip

{\bf Theorem 2.5$\,'$} (the Harnack inequality). {\it Let $\cal L$ be an operator 
of the form (\,{\bf DE}) in $B_{2R}$, and let the condition (\ref{ell}) be satisfied. 
Suppose also that}
\begin{equation}\label{b_large_d}
\|{\bf b}\|_{n,B_{2R}}\le \mathfrak B \quad\mbox{for}\quad n\ge3;
\qquad\qquad \|{\bf b}\ln^{\frac 12}\big(1+R|{\bf b}|\big)\|_{2,B_{2R}}\le \mathfrak B
\quad\mbox{for}\quad n=2.
\end{equation}	
{\it Then there exists a positive constant $N'_3$ depending only on $n$, $\nu$ and $\mathfrak B$ 
such that any nonnegative weak solution of ${\cal L}u=0$ in $B_{2R}$ satisfies}
\begin{equation}\label{Harnack2}
\sup\limits_{B_{R}}u\le N'_3\cdot\inf\limits_{B_{R}}u.
\end{equation}

\begin{proof} We split the spherical layer $B_{2R}\setminus B_{R}$ to $M$ layers of equal 
thickness $\frac R M$ and put $\delta=\frac{1}{2M}$. Obviously, one can choose $M$ depending
only on $n$, $\nu$ and $\mathfrak B$ such that at least for one of these layers (say,
$K=\{r-2\delta R<|x|<r+2\delta R\}$) the following estimates hold:
$$\|{\bf b}\|_{n,K}\le \ep \quad\mbox{for}\quad n\ge3;
\qquad\qquad \|{\bf b}\ln^{\frac 12}\big(1+2\delta R|{\bf b}|\big)\|_{2,K}\le \ep
\quad\mbox{for}\quad n=2$$
(here $\ep=\ep(n,\nu)$ is the above smallness constant).
 
We cover the sphere $|x|=r$ with a finite set of balls $B_{\delta R}(x)$ 
(note that the number of balls depends only on $\delta$). 
Since all doubled balls $B_{2\delta R}(x)$ lye in $K$, we can apply Harnack's inequality in these balls.
This gives $\sup\limits_{|x|=r}u\le C(n,\nu,\delta)\cdot\inf\limits_{|x|=r}u$. However, by the maximum principle,
$$\inf\limits_{B_{R}}u\ge\inf\limits_{|x|=r}u;\qquad\qquad \sup\limits_{B_{R}}u\le\sup\limits_{|x|=r}u,$$ 
and the statement follows.
\end{proof}

The following statement can be proved by verbatim repetition of the proof of Theorem 2.6, using Theorem 2.5$\,'$.\medskip

{\bf Theorem 2.6$\,'$} (the Liouville theorem). {\it Let $\cal L$ be an operator 
of the form (\,{\bf DE}) in $\mathbb R^n$. Let the condition (\ref{ell}) be satisfied,
and}
$${\bf b}\in L_{n,loc}(\mathbb R^n)
\quad\mbox{for}\quad n\ge3;
\qquad\qquad {\bf b}\ln^{\frac 12}\big(1+|{\bf b}|\big)\in L_{2,loc}(\mathbb R^2) 
\quad\mbox{for}\quad n=2.
$$
{\it Suppose also that for some $\delta>0$}
\begin{multline}\label{Liouville2}
\liminf\limits_{R\to\infty}\sup\limits_{|x|=R}
\|{\bf b}\|_{n,B_{\delta R}(x)}<\infty,\quad n\ge3;\\
\liminf\limits_{R\to\infty}\sup\limits_{|x|=R}
\|{\bf b}\ln^{\frac 12}\big(1+R|{\bf b}|\big)\|_{2,B_{\delta R}(x)}<\infty, 
\ \ n=2.
\end{multline}
{\it Then any weak semibounded solution of ${\cal L}u=0$ in $\mathbb R^n$ is a 
constant.}\medskip

{\bf Remark 3}. If ${\bf b}\in L_n(\mathbb R^n)$ (respectively, 
${\bf b}\ln^{\frac 12}\big(1+|x|\,|{\bf b}|\big)\in L_2(\mathbb R^2)$), then (\ref{Liouville2}) is
obviously satisfied.\medskip

As for H\"older estimates for solutions, we have two possibilities. The first one is to take in 
the proof of Lemma 2.3 $R$ sufficiently small, such that the smallness assumptions on ${\bf b}$
are satisfied in $B_{2R}$. This gives the estimate (\ref{Holder}) with $\gamma$ depending only on $n$ 
and $\nu$, while $N_2$ depends also on the moduli of continuity of ${\bf b}$ in $L_n(B_{R_0})$ 
(respectively, of ${\bf b}\ln^{\frac 12}\big(1+R|{\bf b}|\big)$ in $L_2(B_{R_0})$). The second possibility 
is to use Theorem 2.5$\,'$. This gives (\ref{Holder}) with both $\gamma$ and $N_2$ depending on $n$, $\nu$
and $\|{\bf b}\|_{n,B_{R_0}}$ (respectively, 
$\|{\bf b}\ln^{\frac 12}\big(1+R_0|{\bf b}|\big)\|_{2,B_{R_0}}$).

\section{Parabolic case}

{\bf Lemma 3.1}. {\it Let $\cal M$ be an operator 
of the form (\,{\bf DP}) in $Q_{R}^{\lambda,\theta}(x^0;t^0)$, $\lambda>1$, $\theta>0$, and let the conditions (\ref{ell}) and (\ref{bezdiv}) be satisfied. Let also
${\bf b}\in L_{q,\ell}(Q_{R}^{\lambda,\theta}(x^0;t^0))$, with some $q$ and $\ell$ such that
\begin{equation}\label{alpha}
\alpha=\alpha(q,\ell)\equiv\frac nq+\frac 2\ell-1\in\,[0,1[
\end{equation}
($q$ as well as $\ell$ may be infinite\footnote{For $\alpha=0$, the assumption (\ref{bezdiv}) can be removed, with some limitation in the case $q=n$, $\ell=\infty$. We discuss this at the end of this Section.}). 

Then there exists a positive constant $N_4$ depending on $n$, $\nu$, $\lambda$, $\theta$, $q$, $\ell$
and the quantity 
$$\widehat{\cal N}=\widehat{\cal N}(R,\lambda,\theta)\equiv 
R^{-\alpha}\|{\bf b}\|_{q,\ell,Q_{R}^{\lambda,\theta}(x^0;t^0)}$$ 
such that any Lipschitz subsolution of the equation ${\cal M}u=0$ in $Q_{R}^{\lambda,\theta}(x^0;t^0)$ satisfies}
\begin{equation}\label{estmax'}
\sup\limits_{Q_{R}^{1,\frac \theta2}(x^0;t^0)} u_+\le N_4
\bigg(\!\!\!\!\!\!\Xint{\ \ \qquad-}\limits_{Q_{R}^{\lambda,\theta}(x^0;t^0)}\!\!\!u_+^2dxdt\bigg)^{\frac 12}.
\end{equation}
\medskip

{\bf Remark 4}. The quantity $\widehat{\cal N}$ depends also on $q$ and $\ell$. However, we assume these
parameters hold fixed, and we do not indicate this dependence. 
%$$\widehat{\cal N}_{q,\ell}(R,\lambda,\theta)\le (C(n)\lambda^n)^{\frac 1q-\frac 1{\widetilde q}}
%\theta^{\frac 1\ell-\frac 1{\widetilde\ell}}
%\widehat{\cal N}_{\widetilde q,\widetilde\ell}(R,\lambda,\theta).
%$$
We also do not indicate the dependence $\widehat{\cal N}$ on $x^0$ and $t^0$.

\begin{proof} We proceed similarly to Lemma 2.1. Without loss of generality, we assume $(x^0;t^0)=(0;0)$.

For a nonnegative test function $\eta$ we have
$$\int\limits_{Q_{R}^{\lambda,\theta}}(\partial_tu\eta+ a_{ij}D_juD_i\eta+b_iD_iu\,\eta)\,dxdt\le0.
$$
We take $\eta=\varphi'(u)\cdot\xi$, where $\xi$ is a cut-off function, Lipschitz in $x$ and 
vanishing at the neighborhood of $\partial'Q_{R}^{\lambda,\theta}$, while 
$\varphi\in{\cal C}^2(\mathbb R)$ is a convex function vanishing in  $\mathbb R_-$. This gives
\begin{equation}\label{Moser'}
\int\limits_{Q_{R}^{\lambda,\theta}\cap\{u>0\}}\Big(\partial_tv\xi+a_{ij}D_jvD_i\xi+\frac {\varphi''(u)}
{\varphi^{\prime2}(u)}\,a_{ij}D_jvD_iv\,\xi+b_iD_iv\,\xi\Big)\,dxdt\le0,
\end{equation}
where $v=\varphi(u)$.

As in Section 2, by mollification at a neighborhood of the origin, one can weaken in (\ref{Moser}) the assumption
$\varphi\in{\cal C}^2(\mathbb R)$ to $\varphi\in{\cal C}^2(\mathbb R_+\cup\mathbb R_-)$.

Now we put in (\ref{Moser'}) $\varphi(\tau)=\tau_+^p$, $p>1$, and $\xi=\chi_{\{t<\bar t\}}\cdot v\zeta^2$ 
where $\zeta$ is a smooth cut-off function in $Q_{R}^{\lambda,\theta}$, $\bar t \in\,]-\theta R^2,0[$. 
Then we obtain
\begin{multline}\label{Moser_power'}
\frac 12\int\limits_{B_{\lambda R}}(v\zeta)^2\big|^{t=\bar t}dx+\\
+\int\limits_{Q_{R}^{\lambda,\theta}}\chi_{\{t<\bar t\}}\,\Big(\frac {2p-1}{p}\,
a_{ij}D_jvD_iv\,\zeta^2+2a_{ij}D_jv\,vD_i\zeta\,\zeta-v^2\zeta\partial_t\zeta+b_iD_iv\,v\zeta^2\Big)
\,dxdt\le0.
\end{multline}

The last term in (\ref{Moser_power'}) can be estimated using (\ref{bezdiv}) and the H\"older
inequality:
\begin{multline}\label{bezdiv1'}
-\int\limits_{Q_{R}^{\lambda,\theta}}\chi_{\{t<\bar t\}}\,b_iD_iv\,v\zeta^2\,dxdt\le
\int\limits_{Q_{R}^{\lambda,\theta}}\chi_{\{t<\bar t\}}\,b_iv^2\zeta\,D_i\zeta\,dxdt\le\\
\le\|{\bf b}\|_{q,\ell,Q_{R}^{\lambda,\theta}}\|v\zeta\|^{2-\frac 1s}_{r,l,Q_{R}^{\lambda,\theta}}
\|v\zeta^{1-s}|D\zeta|^s\|^{\frac 1s}_{2,2,Q_{R}^{\lambda,\theta}},
\end{multline}
where $s>2$ is defined by $\frac 1s=1-\frac n{2q}-\frac 1\ell$ while $r$ and $l$ are defined by
$$\frac 1{2s}+\frac 1q+\frac {2-\frac 1s}r=1;\qquad \frac 1{2s}+\frac 1\ell+\frac {2-\frac 1s}l=1.
$$
Note that $\frac n2<\frac nr+\frac 2l<\frac n2+1$, and, by the embedding theorem \cite[Ch. II, (3.4)]{LSU},
\begin{equation}\label{embed'}
\|v\zeta\|_{r,l,Q_{R}^{\lambda,\theta}}\le
C_9(n,r,l,\lambda,\theta)R^{\frac nr+\frac 2l-\frac n2}\,\|v\zeta\|_{{\cal V}(Q_{R}^{\lambda,\theta})}.
\end{equation}
Using (\ref{ell}), (\ref{bezdiv1'}) and (\ref{embed'}) and the Young inequality, we obtain from 
(\ref{Moser_power'})
\begin{equation}\label{ee'}
\|v\,\zeta\|^2_{{\cal V}(Q_{R}^{\lambda,\theta})}\le
 C_{10}(n,\nu,q,\ell,\lambda,\theta)\cdot \int\limits_{Q_{R}^{\lambda,\theta}}
v^2\big(|D\zeta|^2+\zeta|\partial_t\zeta|+
R^2\|{\bf b}\|^{2s}_{q,\ell,Q_{R}^{\lambda,\theta}}\zeta^{2-2s}|D\zeta|^{2s}\big)\,dxdt.
\end{equation}

We put $\lambda_m=1+2^{-m}(\lambda-1)$, $\theta_m=\frac \theta2 (1+4^{-m})$, $m\in \mathbb N\cup\{0\}$, and
substitute $\zeta=\zeta_m$ such that 
$$ \zeta_m\equiv1\ \ \mbox{in}\ \ Q_{R}^{\lambda_{m+1},\theta_{m+1}};\quad \zeta_m\equiv0\ \ \mbox{out of}\ \ Q_{R}^{\lambda_m,\theta_m};
\qquad |\partial_t\zeta_m|\le \frac {4^mC}{\theta R^2};\quad \frac{|D\zeta_m|}{\zeta_m^{1-\frac 1s}}\le \frac {2^mC_3(s)}{(\lambda-1)R}.$$ 
Then (\ref{ee'}) implies
\begin{equation}\label{eee'}
\|v\,\zeta_m\|_{{\cal V}(Q_{R}^{\lambda_m,\theta_m})}\le
\frac {C_{11}(n,\nu,q,\ell,\lambda,\theta)}R\cdot\|v\|_{2,2,Q_{R}^{\lambda_m,\theta_m}}\cdot
\big(2^m+\big(2^m\widehat{\cal N}\big)^s\big).
\end{equation}

Now for $p=p_m\equiv (\frac {n+2}n)^m$
 we obtain from (\ref{embed'}) (with $r=l=\frac{2(n+2)}{n}$) and (\ref{eee'})
\begin{multline}\label{iteration'}
\bigg(\!\!\!\!\!\!\!\Xint{\qquad\quad -}\limits_{Q_{R}^{\lambda_{m+1,\theta_{m+1}}}}\!\!\!
u_+^{2p_{m+1}}dxdt\bigg)^{\frac {1}{2p_{m+1}}}\!\le
\bigg(C(n)\!\!\Xint{\quad\ \, -}\limits_{Q_{R}^{\lambda_m,\theta_m}} (v\zeta_m)^r dxdt\bigg)^{\frac {1}{rp_m}}\le\\
\le \bigg(2^{2ms}C_{12}\!\!\Xint{\quad\ \,-}\limits_{Q_{R}^{\lambda_m,\theta_m}}\! v^2\,dxdt\bigg)
^{\frac {1}{2p_m}}\!\!=
\bigg(2^{2ms}C_{12}\!\!\Xint{\quad\ \,-}\limits_{Q_{R}^{\lambda_m,\theta_m}}\! u_+^{2p_m}
dxdt\bigg)^{\frac 1{2p_m}}\!,
\end{multline}
 where $C_{12}$ depends only on $n$, $\nu$, $q$, $\ell$, $\lambda$, $\theta$ and $\widehat{\cal N}$.

Iterating (\ref{iteration'}) we arrive at (\ref{estmax'}).
\end{proof}

{\bf Remark 5}. If ${\bf b}={\bf b}^{(1)}+{\bf b}^{(2)}$, and 
${\bf b}^{(j)}\in L_{q_j,\ell_j}(Q_{R}^{\lambda,\theta}(x^0;t^0))$, $j=1,2$, with 
$\alpha(q_j,\ell_j)\in\,[0,1[$, then the proof of Lemma 3.1
does not change. The same is true for other statements of this Section.\medskip

{\bf Remark 6}. If, under assumptions of Lemma 3.1, $u$ satisfies additionally
\begin{equation}\label{negat}
u(\cdot;t^0-\theta R^2)\le 0\qquad\mbox{in}\quad B_{\lambda R}(x^0),
\end{equation}
then we can estimate $u$ up to the bottom of the cilinder, i.e. one can replace the left-hand side in 
(\ref{estmax'}) by $\sup\limits_{Q_{R}^{1,\theta}(x^0;t^0)} u_+$. Indeed, one may simply put 
$\theta_m\equiv\theta$ and take $\zeta_m$ independent on $t$.\medskip

{\bf Corollary 3.1}. {\it Let $\cal M$ satisfy the assumptions of Lemma 3.1 in $Q_{R}^{\lambda,\theta}(x^0;t^0)$. 

1. If a Lipschitz subsolution of ${\cal M}u=0$ in $Q_{R}^{\lambda,\theta}(x^0;t^0)$ satisfies}
\begin{equation}\label{tiny'}
\mbox{meas}\, (\{u> k\}\cap Q_{R}^{\lambda,\theta}(x^0;t^0))\le \mu\,\mbox{meas}\,(Q_{R}^{\lambda,\theta}),
\qquad \mu< N_4^{-2},
\end{equation}
{\it for some $k$, then}
\begin{equation}\label{estmax1'}
\sup\limits_{Q_{R}^{1,\frac \theta2}(x^0;t^0)} (u-k)\le
N_4\sqrt{\mu}\sup\limits_{Q_{R}^{\lambda,\theta}(x^0;t^0)}(u-k),
\end{equation}
{\it (here $N_4$ is the constant from Lemma 3.1).

2. If a Lipschitz nonnegative supersolution of ${\cal M}V=0$ in $Q_{R}^{\lambda,\theta}(x^0;t^0)$ satisfies}
\begin{equation}\label{tiny1'}
\mbox{meas}\, (\{V< k\}\cap Q_{R}^{\lambda,\theta}(x^0;t^0))\le \mu\,\mbox{meas}\,(Q_{R}^{\lambda,\theta}),
\qquad \mu\le \mu_1\equiv (2N_4)^{-2},
\end{equation}
{\it for some $k>0$, then}
\begin{equation}\label{estmin'}
V\ge \frac k2\qquad\mbox{in}\quad Q_{R}^{1,\frac \theta2}(x^0;t^0).
\end{equation}
{\it If $V$ additionally satisfies}
$$V(\cdot;t^0-\theta R^2)\ge k\qquad\mbox{in}\quad B_{\lambda R}(x^0),
$$
{\it then the estimate (\ref{estmin'}) holds in $Q_{R}^{1,\theta}(x^0;t^0)$.}

\begin{proof} 1. We apply Lemma 3.1 to $u-k$.

2. We apply Lemma 3.1 and Remark 6 to $u=k-V$.
\end{proof}

{\bf Lemma 3.2}. {\it Let $\cal M$ satisfy the assumptions of Lemma 3.1 in $Q_{R}(x^0;t^0)$. 
For any $\delta_0\in\,]0,1]$ there exists $\theta_0\in\,]0,1[$ such that if a Lipschitz nonnegative supersolution of 
${\cal M}V=0$ in $Q_{R}^{1,\theta_0}(x^0;t^0)$ satisfies}
%\begin{equation}
$$\mbox{meas}\, (\{V(\cdot;t^0-\theta_0 R^2)\ge k\}\cap B_{R}(x^0))\ge \delta_0\,\mbox{meas}\,(B_{R})
$$%\end{equation}
{\it for some $k>0$, then}
%\begin{equation}
$$\mbox{meas}\, (\{V(\cdot;\bar t)\ge \frac {\delta_0}3 k\}\cap B_{R}(x^0))\ge
\frac {\delta_0}3 \,\mbox{meas}\,(B_{R})
\qquad\mbox{for any}\quad \bar t\in\,[t^0-\theta_0 R^2,t^0].
$$%\end{equation}
{\it Moreover, $\theta_0$ is completely determined by ${\delta_0}$, $n$, $\nu$, $q$, $\ell$ and the quantity 
$\widehat{\cal N}$.}

\begin{proof}  Without loss of generality, we assume $(x^0;t^0)=(0;0)$. For a nonnegative test function $\eta$ we 
have
\begin{equation}\label{super}
\int\limits_{Q_{R}^{\lambda,\theta_0}}(\partial_tV\eta+ a_{ij}D_jVD_i\eta+b_iD_iV\,\eta)\,dxdt\ge0.
\end{equation}
We take $\eta=\chi_{\{t<\bar t\}}\cdot(V-k)_-\zeta^2(x)$, where $\zeta$ is a smooth cut-off function in $B_{R}$,
$\bar t \in\,]-\theta_0 R^2,0]$. Using (\ref{ell}), (\ref{bezdiv}) and the Young inequality, we obtain
\begin{multline}\label{measure}
\int\limits_{B_{R}}(V-k)_-^2\zeta^2\big|^{t=\bar t}dx
+\nu\int\limits_{Q_{R}^{1,\theta_0}}\chi_{\{t<\bar t\}}\,
|D(V-k)_-|^2\zeta^2dxdt%\le\\
\le \int\limits_{B_{R}}(V-k)_-^2\zeta^2\big|^{t=-\theta_0 R^2}dx+\\
+\int\limits_{Q_{R}^{1,\theta_0}}\chi_{\{t<\bar t\}}\,\Big(C_{13}(n,\nu)(V-k)_-^2|D\zeta|^2dxdt+
2b_i(V-k)_-^2\,\zeta\,D_i\zeta\Big)dxdt.
\end{multline}
Now we choose $\zeta$ such that $\zeta\equiv1$ in $B_{(1-\sigma)R}$ and $|D\zeta|\le \frac {2}{\sigma R}$ where $\sigma<1$ is a parameter to be chosen later. Observing that $(V-k)_-^2\le k^2$, we estimate the right-hand side of (\ref{measure}) by
$$k^2\Big[(1-{\delta_0})\,\mbox{meas}\,(B_R)+C_{13}\theta_0 R^2\cdot \frac{4\mbox{meas}\,(B_R)}{(\sigma R)^2}+
\frac{2}{\sigma R}\,\|{\bf b}\|_{q,\ell,Q_{R}}\|{\bf 1}\|_{q',\ell',Q_{R}^{1,\theta_0}}\Big].
$$
On the another hand, 
$$\int\limits_{B_{R}}(V-k)_-^2\zeta^2\big|^{t=\bar t}dx\ge\!\!
\int\limits_{\{V<\frac {\delta_0}3 k\}\cap B_{(1-\sigma)R}}\!\!\!\!(V-k)_-^2\big|^{t=\bar t}dx\ge
\big(1-\frac {\delta_0}3\big)^2 k^2\, \mbox{meas}\, (\{V(\cdot;\bar t)< \frac {\delta_0}3 k\}\cap B_{(1-\sigma)R}).
$$
Thus,
$$\mbox{meas}\, (\{V(\cdot;\bar t)< \frac {\delta_0}3 k\}\cap B_{(1-\sigma)R})\le 
\frac {\mbox{meas}\,(B_R)}{(1-\frac {\delta_0}3)^2}\cdot\Big[(1-{\delta_0})+\frac{4C_{13}\theta_0}{\sigma^2}+
\frac{C(n)\theta_0^{\frac 2{\ell'}}{\widehat{\cal N}}}{\sigma}\,\Big],
$$
and therefore, 
$$\mbox{meas}\, (\{V(\cdot;\bar t)< \frac {\delta_0}3 k\}\cap B_{R})\le 
\frac {\mbox{meas}\,(B_R)}{(1-\frac {\delta_0}3)^2}\cdot
\Big[(1-{\delta_0})+C(n)\sigma+\frac{4C_{13}\theta_0}{\sigma^2}+
\frac{C(n)\theta_0^{\frac 2{\ell'}}{\widehat{\cal N}}}{\sigma}\,\Big].
$$
Since $1-{\delta_0}\le(1-\frac {\delta_0}3)^3-\frac 8{27}{\delta_0}^2$, one can choose $\sigma$ and then $\theta_0$ small enough such that the right-hand side is not greater that $(1-\frac {\delta_0}3)\, \mbox{meas}\,(B_R)$, and the Lemma follows.
\end{proof}

{\bf Lemma 3.3}. {\it Let $\cal M$ satisfy the assumptions of Lemma 3.1 in $Q_{R}^{\lambda,\theta}(x^0;t^0)$ with $\lambda>1$. Let a Lipschitz nonnegative supersolution of ${\cal M}V=0$ in $Q_{R}^{\lambda,\theta}(x^0;t^0)$ satisfy}
\begin{equation}\label{measure1}
\mbox{meas}\, (\{V(\cdot;t)\ge k_0\}\cap B_{R}(x^0))\ge \delta_1\,\mbox{meas}\,(B_{R})
\qquad\mbox{for any}\quad t\in\,[t^0-\theta R^2,t^0]
\end{equation}
{\it for some $k_0>0$ and $\delta_1>0$. Then for any $\mu\in\,]0,1[$ there exists $s>1$ such that}
$$\mbox{meas}\, (\{V< 2^{-s} k_0\}\cap Q_{R}^{1,\theta}(x^0;t^0))\le \mu\,\mbox{meas}\,(Q_{R}^{1,\theta}).
$$
{\it Moreover, $s$ is completely determined by $n$, $\nu$, $\lambda$, $\theta$, $\mu$, $\delta_1$, $q$,
$\ell$, and the quantity $\widehat{\cal N}$.}

\begin{proof}  Without loss of generality, we assume $(x^0;t^0)=(0;0)$. For $m\in{\mathbb Z}_+$ we put $k_m=2^{-m}k_0$,
$${\cal E}_m(t)=\{x\in B_R:\ k_{m+1}\le V(x,t)< k_m\};\qquad
{\cal E}_m=\{(x;t):\ t\in\,[-\theta R^2,0],\ x\in{\cal E}_m(t)\}.
$$
We take in (\ref{super}) $\eta=(V-k_m)_-\zeta^2(x)$, where $\zeta$ is a smooth cut-off function, vanishing at the neighborhood of $\partial B_{\lambda R}$ and satisfying $\zeta\equiv 1$ in $B_R$, 
$|D\zeta|\le \frac 2{(\lambda-1)R}$. Similarly to the proof of Lemma 3.2, we derive
\begin{equation}\label{fff}
\int\limits_{\{V<k_m\}}
|DV|^2\zeta^2dxdt
=\int\limits_{Q_{R}^{\lambda,\theta}}
|D(V-k_m)_-|^2\zeta^2dxdt\le C_{14}(n,\nu,\lambda, \theta,\ell,\widehat{\cal N}) k_m^2 R^n.
\end{equation}

Further, De Giorgi's inequality (see, e.g., \cite[Ch. II, (5.6)]{LSU}) and the assumption (\ref{measure1}) give
$$(k_m-k_{m+1})\cdot \mbox{meas}\, (\{V(\cdot;t)< k_{m+1}\}\cap B_{R})\le 
\frac{C(n)R}{\delta_1}\int\limits_{{\cal E}_m(t)}|DV(\cdot;t)|\,dx,\qquad t\in\,[-\theta R^2,0].
$$
We integrate this relation w.r.t $t$ and then square both parts, arriving at
$$k_{m+1}^2\,\mbox{meas}^2 (\{V< k_{m+1}\}\cap Q_{R}^{1,\theta})\le 
\frac{C(n)R^2}{\delta_1^2}\int\limits_{{\cal E}_m}|DV|^2 dxdt\cdot\mbox{meas}\,({\cal E}_m).
$$
Together with (\ref{fff}), this gives
$$\mbox{meas}^2 (\{V< k_{m+1}\}\cap Q_{R}^{1,\theta})\le 
C(n)C_{14}\delta_1^{-2}R^{n+2}\cdot\mbox{meas}\,({\cal E}_m).
$$
Therefore,
\begin{multline*}
s\cdot\mbox{meas}^2 (\{V< k_s\}\cap Q_{R}^{1,\theta})\le 
\sum\limits_{m=0}^{s-1}\mbox{meas}^2 (\{V< k_{m+1}\}\cap Q_{R}^{1,\theta})\le \\
\le C_{15}\delta_1^{-2}\cdot\mbox{meas}\,(Q_{R}^{1,\theta})\cdot\sum\limits_{m=0}^{s-1}\mbox{meas}\,({\cal E}_m)\le
C_{15}\delta_1^{-2}\cdot\mbox{meas}^2(Q_{R}^{1,\theta})
\end{multline*}
(here $C_{15}$ depends on the same quantities as $C_{14}$), and the Lemma follows.
\end{proof}

{\bf Corollary 3.2}. {\it Let $\cal M$ satisfy the assumptions of Lemma 3.1 in $Q_{R}^{2,1}$. 

1. If a Lipschitz nonnegative supersolution of ${\cal M}V=0$ in $Q_{R}^{2,1}$ satisfies}
\begin{equation}\label{measure2}
\mbox{meas}\, (\{V(\cdot;\bar t)\ge k\}\cap B_{R})\ge \delta\,\mbox{meas}\,(B_{R})
\end{equation}
{\it for some $\bar t\in\,[-R^2, -\Theta R^2]$, $k>0$, $\delta,\Theta\in\,]0,1]$, then}
\begin{equation}\label{estmin1'}
V\ge \beta_1 k\qquad\mbox{in}\quad Q_{R}^{1,\frac {\theta_1}2}(0;\bar t+\theta_1 R^2).
\end{equation}
{\it Here $\theta_1=\min\{\Theta,\theta_0\}$, where 
$\theta_0=\theta_0(\delta,n,\nu,q,\ell,\widehat{\cal N})$ is the constant from Lemma 3.2,
while $\beta_1$ depends on the same quantities as $\theta_0$.

2. If the relation (\ref{measure2}) holds with $\delta=1$, i.e.}
\begin{equation}\label{estmin2'}
V(\cdot;\bar t)\ge k\qquad\mbox{in}\quad B_{R},
\end{equation}
{\it then for any $\sigma\in\,]0,1[$ 
\begin{equation}\label{estmin1''}
V\ge\beta_1 k\qquad\mbox{in}\quad Q_{R}^{\sigma,\theta_1}(0;\bar t+\theta_1 R^2),
\end{equation}
In this case $\beta_1$ depends additionally on $\sigma$.}

\begin{proof} First, we use Lemma 3.2 with $t^0=\bar t+\theta_1 R^2$. Then, in the case 1, we apply Lemma
3.3 with $R\to\frac 32 R$, $\lambda=\frac 43$, $\theta=\frac 49\theta_1$,
$\delta_1=\frac{2^n\delta}{3^{n+1}}$, $k_0=\frac{\delta}3 k$ and $\mu=\mu_1$, where 
$\mu_1=\mu_1(n,\nu,\frac 32, \theta_1,q,\ell,\widehat{\cal N})$ is the constant from Corollary 3.1, 
part 2. Finally, Corollary 3.1 with $\lambda=\frac 32$ and  $\theta=\theta_1$ 
gives (\ref{estmin1'}) with $\beta_1=\frac {\delta}{3\cdot 2^{s+1}}$, where 
$s=s(n,\nu,\frac 32, \theta_1, \mu_1,\delta_1,q,\ell,\widehat{\cal N})$ is the constant from Lemma 3.3.

In the case 2, we apply Lemma 3.3 with $\lambda=2$, $\theta=\theta_1$, $\delta_1=\frac{\delta}3$, 
$k_0=\frac{\delta}3 k$ and $\mu=\mu_1(n,\nu,\sigma^{-1}, \sigma^{-2}\theta_1,q,\ell,\widehat{\cal N})$. Finally, the last statement of Corollary 3.1 with $R\to\sigma R$, $\lambda=\sigma^{-1}$ and $\theta=\sigma^{-2}\theta_1$ gives (\ref{estmin1''}) with $\beta_1=\frac {\delta}{3\cdot 2^{s+1}}$, where 
$s=s(n,\nu,\sigma^{-1}, \sigma^{-2}\theta_1, \mu_1,\delta_1,q,\ell,\widehat{\cal N})$.
\end{proof}

{\bf Lemma 3.4}. {\it Let $\cal M$ satisfy the assumptions of Lemma 3.1 in $Q_{R}^{2,1}$. Let a Lipschitz nonnegative supersolution of ${\cal M}V=0$ in $Q_{R}^{2,1}$ satisfy (\ref{estmin2'}) for some $k>0$ and 
$\bar t\in\,[-R^2,-\Theta R^2]$, $\Theta\in\,]0,1]$. Then}
\begin{equation}\label{estmin3'}
V\ge \beta_2 k\qquad\mbox{in}\quad \widehat Q=B_{\frac R2}\times[\bar t,0].
\end{equation}
{\it Moreover, $\beta_2$ is completely determined by $\Theta$, $n$, $\nu$, $q$, $\ell$, 
and the quantity $\widehat{\cal N}$.}

\begin{proof} 
We set $M=\mbox{entier}\big(\frac{|\bar t|}{\theta_1 R^2}\big)+1$ and 
$\widehat\theta_1=\frac {|\bar t|}{MR^2}$. 
Now let us consider cylinders 
$$Q^{(m)}=Q_{R}^{1-\frac m{2M},\widehat\theta_1}(0;\bar t+m\widehat\theta_1R^2), \qquad m=1,\dots,M.$$ 
By (\ref{estmin1''}) we consequently obtain 
$$V\ge\widehat\beta_1\cdot\inf\limits_{Q^{(m)}} V \qquad\mbox{in}\quad Q^{(m+1)},
$$
where $\widehat\beta_1=\widehat\beta_1(n,\Theta,\nu,q,\ell,\widehat{\cal N})>0$.

Since $\widehat Q\subset\bigcup\limits_m Q^{(m)}$, this ensures (\ref{estmin3'}) with $\beta_2=\widehat\beta_1^M$.
\end{proof}

{\bf Corollary 3.3}. {\it Let $\cal M$ satisfy the assumptions of Lemma 3.1 in $Q_{R}^{2,1}$. 
If a Lipschitz nonnegative supersolution of ${\cal M}V=0$ in $Q_{R}^{2,1}$ satisfies (\ref{measure2})
for some $k>0$, $\delta\in\,]0,1]$ and $\bar t\in\,[-R^2, -\frac 34 R^2]$, then}
\begin{equation}\label{estmin4'}
V\ge \beta_3 k\qquad\mbox{in}\quad Q_{\frac R2},
\end{equation}
{\it where $\beta_3$ is completely determined by $\delta$, $n$, $\nu$, $q$, $\ell$, 
and the quantity $\widehat{\cal N}$.}

\begin{proof} It suffices to apply consequently Corollary 3.2, part 1, and Lemma 3.4. 
\end{proof}

{\bf Corollary 3.4}. {\it Let $\cal M$ satisfy the assumptions of Lemma 3.1 in $Q_{R}^{2,1}$. 
If a Lipschitz nonnegative supersolution of ${\cal M}V=0$ in $Q_{R}^{2,1}$ satisfies}
\begin{equation}\label{measure3}
\mbox{meas}\, (\{V>k\}\cap Q_{R})\ge \widehat\delta\,\mbox{meas}\,(Q_{R}),
\end{equation}
{\it for some $k>0$, $\widehat\delta\in\,]0,1]$, then}
\begin{equation}\label{estmin5'}
V\ge \beta_4 k\qquad\mbox{in}\quad Q_{R}^{\frac 12,\frac {\widehat\delta}4},
\end{equation}
{\it where $\beta_4$ is completely determined by $\widehat\delta$, $n$, $\nu$, $q$, $\ell$, 
and the quantity $\widehat{\cal N}$.}
 
\begin{proof} The inequality (\ref{measure3}) obviously implies
$$\mbox{meas}\, (\{V>k\}\cap Q_{R}^{1,1-\frac{\widehat\delta}2}(0;-\frac{\widehat\delta}2R^2))\ge 
\frac {\widehat\delta}2\,\mbox{meas}\,(Q_{R}).$$ 
Therefore, there exists 
$\bar t\in\,[-R^2,-\frac{\widehat\delta}2R^2]$, such that (\ref{measure2}) holds with 
$\delta=\frac{\widehat\delta}2$. By Corollary 3.2, part 1, 
$V(\cdot;\bar t+\frac{\theta_1}2R^2)\ge \beta_1 k$ in $B_{R}$. Finally, we observe that 
$\bar t+\frac{\theta_1}2R^2\le-\frac{\widehat\delta}4R^2$, and Lemma 3.4 provides (\ref{estmin5'}).
\end{proof}

{\bf Corollary 3.5} (strong maximum principle). {\it Let $\cal M$ satisfy the assumption of Lemma 3.1
in $Q$. Then any Lipschitz nonconstant supersolution of ${\cal M}V=0$ in $Q$ cannot attain its minimum 
at a point of $\partial Q\setminus\partial'Q$.}\medskip

\begin{proof} Without loss of generality, $\inf\limits_Q V=0$.

Assume the converse. Then there exists $(x^0;t^0)\in\overline Q\setminus\partial'Q$ such that $V(x^0;t^0)=0$ 
but $V\not\equiv0$ in $Q_R(x^0;t^0)\subset Q_{R}^{2,1}(x^0;t^0)\subset Q$ with some $R$. Then the relation
(\ref{measure3}) holds for some $k>0$ and $\delta>0$, and we obtain (\ref{estmin5'}), a contradiction.
\end{proof}

{\bf Lemma 3.5}. {\it Let $\cal M$ satisfy the assumptions of Lemma 3.1 in $Q_{2R}$. 
Then any Lipschitz solution of ${\cal M}u=0$ in $Q_{2R}$ satisfies the estimate}
\begin{equation}\label{osc1}
\underset{Q_{\frac R2}}{\mbox{osc}}\ u\le \varkappa_1\,\underset{Q_{2R}}{\mbox{osc}}\ u,
\end{equation}
{\it where $\varkappa_1<1$ depends on $n$, $\nu$, $q$ $\ell$ and the quantity $\widehat{\cal N}$.}

\begin{proof} We set $k=\frac 12\underset{Q_{2R}}{\mbox{osc}}\ u$ and consider two functions $V_1=u-\inf\limits_{Q_{2R}}u$ and $V_2=\sup\limits_{Q_{2R}}u-u$. At least one of them satisfies (\ref{measure2})
with $\delta=\frac 12$ and $\bar t=-R^2$. Therefore, Corollary 3.3 gives for this function the estimate
(\ref{estmin4'}), which implies (\ref{osc1}) with $\varkappa_1=1-\frac 12\beta_3(\frac 12, n,\nu,q,\ell,\widehat{\cal N})$.
\end{proof}

{\bf Corollary 3.6} (H\"older estimate). {\it Let $\cal M$ satisfy the assumption of Lemma 3.1
in $Q_{R_0}$. Let also $\sup\limits_{R<R_0}{\widehat{\cal N}}(R,1,1)<\infty$. Then any Lipschitz solution 
of ${\cal M}u=0$ in $Q_{R_0}$ satisfies the estimate}
\begin{equation}\label{Holder'}
 \underset{Q_\rho}{\mbox{osc}}\ u\le N_5 \Big(\frac \rho r\Big)^{\gamma_1}
\cdot\underset{Q_r}{\mbox{osc}}\ u,\qquad 0<\rho<r\le R_0,
\end{equation}
{\it where $N_5$ and $\gamma_1$ depend on $n$, $\nu$, $q$ and 
$\sup\limits_{R<R_0}{\widehat{\cal N}}(R,1,1)$.}

\begin{proof} Iteration of (\ref{osc1}) gives (\ref{Holder'}) with $\gamma_1=-\log_4(\varkappa_1)$.
\end{proof}

{\bf Corollary 3.7} (the Liouville theorem). {\it Let $\cal M$ be an operator 
of the form (\,{\bf DP}) in $\mathbb R^n\times \mathbb R_-$, and let the conditions (\ref{ell}) and 
(\ref{bezdiv}) be satisfied. Let also ${\bf b}\in L_{q,\ell,loc}(\mathbb R^n\times\mathbb R_-)$, with some 
$q$ and $\ell$ satisfying (\ref{alpha}) ($q$ as well as $\ell$ may be infinite). Finally, assume that}
\begin{equation}\label{Liouville'}
\liminf\limits_{R\to\infty}{\widehat{\cal N}}(R,1,1)<\infty.
\end{equation}
{\it Then any Lipschitz bounded solution of ${\cal M}u=0$ in $\mathbb R^n\times\mathbb R_-$ is a 
constant.}\medskip

{\bf Remark 7}. If ${\bf b}\in L_{q,\ell}(\mathbb R^n\times\mathbb R_-)$, then (\ref{Liouville'}) is 
obviously satisfied.

\begin{proof} Iteration of (\ref{osc1}) with respect to a suitable sequence $R_m\to\infty$ gives 
the statement.
\end{proof}

To prove the Harnack inequality, we need the following modification of Lemma 3.4.\medskip

{\bf Lemma 3.4$\,'$} (slant cylinder). {\it Let $\cal M$ satisfy the assumptions of Lemma 3.1 in 
$Q_{R}^{\lambda,1}$, $\lambda>2$. Let a Lipschitz nonnegative supersolution of ${\cal M}V=0$ in 
$Q_{R}^{\lambda,1}$ satisfy} 
$$V(\cdot;-R^2)\ge k\qquad\mbox{in}\quad B_{R}(x^0),
$$
{\it for some $k>0$ and $x^0\in B_{(\lambda-2)R}$. Then for any $x^1\in B_{(\lambda-2)R}$}
\begin{equation}\label{estmin3''}
V\ge \widehat\beta_2 k\qquad\mbox{in}\quad 
\widetilde Q=\{(x;t):\ t\in\,[-R^2,0],\ x\in B_{\frac R2}\big(x^1+(x^1-x^0)\frac t{R^2}\big)\}.
\end{equation}
{\it Moreover, $\widehat\beta_2$ is completely determined by $\lambda$, $n$, $\nu$, $q$, $\ell$, 
and the quantity $\widehat{\cal N}$.}

\begin{proof} We put $\widehat x(t)=x^1+(x^1-x^0)\frac t{R^2}$. Then it is easy to see that the 
function $\widetilde V(x;t)=V(x-\widehat x(t);t)$ is a Lipschitz nonnegative supersolution of 
$$\widetilde{\cal M}V\equiv \partial_tV-D_i\big(\widetilde a_{ij}(x;t)D_jV\big)+\widetilde b_i(x;t)D_iV=0
$$ 
in $Q_{R}^{2,1}$, where
$$\widetilde a_{ij}(x;t)=a_{ij}(x-\widehat x(t);t);\qquad
\widetilde b_i(x;t)=b_i(x-\widehat x(t);t)+\frac{x^1_i-x^0_i}{R^2}.
$$
Note that $\widetilde{\cal M}$ satisfies the assumptions of Lemma 3.1 in $Q_{R}^{2,1}$, and the quantity
$\widetilde{\widehat{\cal N}}$ is bounded by $\widehat{\cal N}+C_{16}(\lambda,n)$. By Lemma 3.4 (with $\Theta=1$), we obtain 
(\ref{estmin3''}).
\end{proof}

The next statement is a parabolic analog of Lemma 2.4. For $(x;t)\in Q$ we introduce the notation 
$d_{\mbox{par}}((x;t),\partial'Q)=\inf \{\rho>0:\ Q_\rho(x;t)\subset Q\}$.\medskip

{\bf Lemma 3.6}. {\it Let $\cal M$ satisfy the assumptions of Lemma 3.1 in $Q_{2R}$, and let 
${\bf b}\in\mathbb M^{\alpha}_{q,\ell}(Q_{2R})$. Let for a Lipschitz nonnegative supersolution of 
${\cal M}V=0$ in $Q_{2R}$ and some $(y;s)\in Q_{R}^{2,2}(0;-2R^2)$, the inequality 
$\inf\limits_{B_\rho(y)}V(\cdot;s)=k>0$ holds with $\rho=\frac 14 d_{\mbox{\rm par}}((y;s),\partial' Q_{2R})$.
Then}
\begin{equation}\label{estmin6'}
\inf\limits_{Q_{R}} V\ge N_6\Big(\frac \rho R\Big)^{\widehat\gamma_1} k,
\end{equation}
{\it where $N_6$ and $\widehat\gamma_1$ depend on $n$, $\nu$, $q$, $\ell$ and 
$\|{\bf b}\|_{\mathbb M^{\alpha}_{q,\ell}(Q_{2R})}$.}

\begin{proof} We denote by ${\mathfrak N}$ an integer number such that 
$2^{-({\mathfrak N}+1)}R\le\rho< 2^{-{\mathfrak N}}R$ and introduce a sequence of cylinders 
$Q_{{\mathfrak r}_m}^{4,1}(y^m;t^m)$, $m=0,\dots,{\mathfrak N}$, as follows:
$$\begin{array}{lll}
{\mathfrak r}_0=2^{-({\mathfrak N}+1)}R,&y^0=y,& t^0=s+{\mathfrak r}_0^2;\\
{\mathfrak r}_m=2{\mathfrak r}_{m-1},&y^m=y^{m-1}-\min\{2{\mathfrak r}_m,|y^{m-1}|\}\,{\bf e},&
t^m=t^{m-1}+{\mathfrak r}_m^2
\end{array}
$$
 (here ${\bf e}=\frac y{|y|}$). Also we denote $y^{-1}=y$, $t^{-1}=s$.

Direct computation shows that $Q_{{\mathfrak r}_m}^{4,1}(y^m;t^m)\subset Q_{2R}$ for all $m=0,\dots,{\mathfrak N}$.
Therefore, the assumptions of Lemma 3.4$\,'$ (with $\lambda=4$, $x^0=y^{m-1}$) are fulfilled in 
$Q_{{\mathfrak r}_m}^{4,1}(y^m;t^m)$. Using Lemma 3.4$\,'$ with $x^1\in B_{2{\mathfrak r}_m}(y^m)$, we obtain the inequality
$$V(\cdot;t^m)\ge \widehat\beta_2(4,n,\nu,q,\ell, \|{\bf b}\|_{\mathbb M^{\alpha}_{q,\ell}(Q_{2R})})
\cdot \inf\limits_{B_{{\mathfrak r}_m}(y^{m-1})} V(\cdot;t^{m-1})
\qquad\mbox{in}\quad B_{\frac 52{\mathfrak r}_m}(y^m)$$
and, in particular, 
$$\inf\limits_{B_{{\mathfrak r}_{m+1}}(y^m)} V(\cdot;t^m)\ge 
\widehat\beta_2\cdot\inf\limits_{B_{{\mathfrak r}_m}(y^{m-1})} V(\cdot;t^{m-1}).$$

Since ${\mathfrak r}_{{\mathfrak N}+1}=R$ and $y^{\mathfrak N}=0$, we obtain
$$\inf\limits_{B_{R}} V(\cdot;t^{\mathfrak N})\ge 
\widehat\beta_2^{{\mathfrak N}+1}\cdot\inf\limits_{B_{{\mathfrak r}_0}(y)} V(\cdot;s)\ge 
\Big(\frac \rho {2R}\Big)^{\widehat\gamma_1}\cdot k,$$
where $\widehat\gamma_1=-\log_2(\widehat\beta_2)$.

It is easy to estimate $t^{\mathfrak N}$:
$$t^{\mathfrak N}=s+\sum\limits_{m=0}^{\mathfrak N}{\mathfrak r}_m^2=
s+{\mathfrak r}_0^2\cdot \frac{2^{2{\mathfrak N}+2}-1}3<s+\frac{R^2}3\le -\frac 53 R^2.
$$

Now we use Lemma 3.4$\,'$ (with $\lambda=4$, $x^0\in B_{\frac R2}$, $x^1\in B_R$) in 
$Q_{\frac R2}^{4,1}(0;t^{\mathfrak N}+\frac{R^2}4)$. Since slant cylinders for $x^0\in B_{\frac R2}$, $x^1\in B_R$
cover $Q_{R}^{1,\frac 18}(0;t^{\mathfrak N}+\frac{R^2}4)$,
we obtain
$$\inf\limits_{Q_{R}^{1,\frac 18}}(0;t^{\mathfrak N}+\frac{R^2}4) V\ge 
\widehat\beta_2\Big(\frac \rho {2R}\Big)^{\widehat\gamma_1}\cdot k.
$$
Repeating this process, we cover $Q_{R}$. This gives (\ref{estmin6'}) with
$N_6=2^{-\widehat\gamma_1}\widehat\beta_2^{31}$.
\end{proof}

{\bf Theorem 3.7} (the Harnack inequality). {\it Let $\cal M$ be an operator 
of the form (\,{\bf DP}) in $Q_{2R}$, and let the conditions (\ref{ell}) and 
(\ref{bezdiv}) be satisfied. Let also ${\bf b}\in \mathbb M^{\alpha}_{q,\ell}(Q_{2R})$, with some 
$q$ and $\ell$ satisfying (\ref{alpha}) (as $q$ as $\ell$ may be infinite). 

Then there exists a positive constant $N_7$ depending on $n$, $\nu$, $q$, $\ell$
and $\|{\bf b}\|_{\mathbb M^{\alpha}_{q,\ell}(Q_{2R})}$, such that any Lipschitz nonnegative solution 
of ${\cal M}u=0$ in $Q_{2R}$ satisfies}
\begin{equation}\label{Harnack'}
\sup\limits_{Q_{R}(0;-2R^2)}u\le N_7\cdot\inf\limits_{Q_{R}}u.
\end{equation}

\begin{proof} Similarly to Theorem 2.5, we denote by $(y;s)$ a maximum point of the 
function
$$v(x;t)=(d_{\mbox{par}}((x;t),\partial'Q_{2R}))^{\widehat\gamma_1}\!\cdot u(x;t); 
\qquad (x;t)\in Q_{R}^{2,2}(0;-2R^2)$$
(here $\widehat\gamma_1$ is the constant from Lemma 3.6) and set
$$\rho=\frac 14\,d_{\mbox{par}}((y;s),\partial'Q_{2R});\qquad 
{\mathfrak M}=v(y;s)=(4\rho)^{\widehat\gamma_1}\!\cdot u(y;s).
$$

It is obvious that
\begin{eqnarray}
\sup\limits_{Q_{R}(0;-2R^2)}u\le \frac {\mathfrak M}{R^{\widehat\gamma_1}}=
\Big(\frac {4\rho} R\Big)^{\widehat\gamma_1}\!\cdot u(y;s);\label{sup'}\\
\sup\limits_{Q_{2\rho}(y;s)}u\le \frac {\mathfrak M}{(2\rho)^{\widehat\gamma_1}}=
2^{\widehat\gamma_1}\cdot u(y;s).\phantom{u(y)}\label{sup1'}
\end{eqnarray}

Denote $k=\frac 12 u(y;s)$. If
$\mbox{meas}\, (\{u> k\}\cap Q_{2\rho}(y;s))\le \mu\,\mbox{meas}\,(Q_{2\rho})$,
then Corollary 3.1, part 1 (with $\lambda=2$ and $\theta=4$) and (\ref{sup1'}) imply the relation
$$k=u(y;s)-k\le \sup\limits_{Q_\rho(y;s)}(u-k)\le N_4\sqrt{\mu}\sup\limits_{Q_{2\rho}(y;s)}(u-k)
\le N_4\sqrt{\mu}\,(2^{\widehat\gamma_1+1}-1)k,
$$
which is impossible for $\mu\le\mu_2\equiv\frac 1{2^{2\widehat\gamma_1+2}}\,N_4^{-2}$. Thus,
$\mbox{meas}\, (\{u> k\}\cap Q_{2\rho}(y;s))\ge \mu_2\,\mbox{meas}\,(Q_{2\rho})$, and Corollary 3.4
(with $\widehat\delta=\mu_2$) gives
\begin{equation}\label{inf'}
\inf\limits_{B_\rho(y)}u(\cdot;s)\ge \beta_4k=\frac{\beta_4}2\cdot u(y;s).
\end{equation}

Finally, Lemma 3.6 gives
\begin{equation}\label{inf1'}
\inf\limits_{Q_{R}}u\ge N_6\Big(\frac \rho R\Big)^{\widehat\gamma_1}\inf\limits_{B_\rho(y)}u(\cdot;s).
\end{equation}
Combining (\ref{sup'}), (\ref{inf'}) and (\ref{inf1'}), we arrive at (\ref{Harnack'}) with 
$N_7=\frac{2^{2\widehat\gamma_1+1}}{\beta_4N_6}$.
\end{proof}

As in Section 2, the proofs of Lemmas 3.1--3.3 run without changes also for weak sub/supersolutions of ${\cal M}u=0$ 
if the bilinear form 
$$\widehat{\cal B}\big\langle u,\eta\big\rangle\equiv \int\limits_{Q_{R}^{\lambda,\theta}}b_iD_iu\,\eta\,dxdt$$
can be continuously extended to the pair $(v,v\zeta^2)$ with $v\in {\cal V}(Q_{R}^{\lambda,\theta})$. 
Unfortunately, we have no parabolic analog of sharp results by Maz'ya--Verbitsky, so we can give only rather rough sufficient conditions. The simplest one is
$${\bf b}\in \mathbb M^{\frac nq-1}_{q,\infty}(Q_{R}^{\lambda,\theta}),\quad \frac n2<q\le n,\quad
\mbox{div}({\bf b})=0.
$$
If Lemmas 3.1--3.3 are proved, all subsequent statements obviously hold true.\medskip

In addition, let us consider the case $\alpha(q,\ell)=0$, i.e. $\frac nq+\frac 2\ell=1$. As in the elliptic case, 
main results of this section hold true for weak (sub/super)solutions without the assumption
(\ref{bezdiv})\footnote{Also the assumptions on ${\bf b}$ can be weakened in the scale of Lorentz spaces.}. 
The only exceptional situation is $q=n$, where the assumption (\ref{bezdiv}) seems to be unavoidable without
the smallness restriction on ${\bf b}$ \footnote{By Remark 4, if the assumption (\ref{bezdiv}) holds, the case 
$q=n$, $\ell=\infty$ is simply included into the case $q=n$, $2<\ell<\infty$.}.  We explain briefly
changes in the proofs.\medskip

Similarly to Lemma 2.1, Lemma 3.1 in this case can be proved under additional assumption
$\|{\bf b}\|_{q,\ell,Q_{R}^{\lambda,\theta}}\le\ep(n,\nu)$.

Lemmas 3.2 and 3.3 are proved without changes. Therefore, all subsequent statements hold true under assumption
of sufficient smallness of  ${\bf b}$.\medskip

In what follows we will assume $q>n$. Then, as in the elliptic case, strong maximum principle 
holds without smallness assumption on ${\bf b}$.

When proving the Harnack inequality, one can exclude the smallness assumption on ${\bf b}$ similarly to Theorem
2.5$\,'$. The result reads as follows.\medskip

{\bf Theorem 3.7$\,'$} (the Harnack inequality). {\it Let $\cal M$ be an operator 
of the form (\,{\bf DE}) in $Q_{2R}$, and let the condition (\ref{ell}) be satisfied. 
Suppose also that}
$$\|{\bf b}\|_{q,\ell,Q_{2R}}\le \mathfrak B,\qquad \frac nq+\frac 2\ell=1,\quad q>n.
$$
{\it Then there exists a positive constant $N'_7$ depending only on $n$, $\nu$, $q$ and $\mathfrak B$ 
such that any nonnegative weak solution of ${\cal M}u=0$ in $Q_{2R}$ satisfies}
$$\sup\limits_{Q_{R}(0;-2R^2)}u\le N'_7\cdot\inf\limits_{Q_{R}}u.
$$

As in the elliptic case, we have two possibilities to prove the H\"older estimates for solutions. 
The first one is to take in Lemma 3.5 $R$ sufficiently small, such that the smallness assumptions 
on ${\bf b}$ are satisfied in $Q_{2R}$. This gives the estimate (\ref{Holder'}) with $\gamma_1$ depending 
only on $n$, $\nu$ and $q$, while $N_5$ depends also on the moduli of integral continuity of ${\bf b}$ in $L_{q,\ell}(Q_{R_0})$. The second possibility is to use Theorem 3.7$\,'$. This gives (\ref{Holder'}) with both 
$\gamma_1$ and $N_5$ depending on $n$, $\nu$, $q$ and $\|{\bf b}\|_{q,\ell,Q_{R_0}}$.\medskip

Finally, the next statement directly follows from the second variant of the H\"older estimate.\medskip

{\bf Corollary 3.7$\,'$} (the Liouville theorem). {\it Let $\cal M$ be an operator of the form (\,{\bf DP}) 
in $\mathbb R^n\times \mathbb R_-$, and let the conditions (\ref{ell}) be satisfied. Let also 
${\bf b}\in L_{q,\ell}(\mathbb R^n\times\mathbb R_-)$, with some $q$ and $\ell$ such that 
$\frac nq+\frac 2\ell=1$, $q>n$.
Then any weak bounded solution of ${\cal M}u=0$ in $\mathbb R^n\times\mathbb R_-$ is a constant.}

\section{Application to a problem of hydrodynamics}

When considering axisymmetric flows of viscous incompressible liquid, the following equation of ({\bf DP}) form arises:
\begin{equation}\label{NS1}
{\cal M}u\equiv \partial_tu-\Delta u+b_i(x',x_3;t)D_iu=0 \qquad\mbox{in}\quad \mathbb R^3\times \mathbb R_-.
\end{equation}
Here we denote $x'=(x_1,x_2)$; 
\begin{equation}\label{NS2}
{\bf b}={\bf v}+\widehat{\bf b}=\Big(v^1+\ep\frac{2x_1}{|x'|^2},v^2+\ep\frac{2x_2}{|x'|^2},v^3\Big),
\end{equation}
where ${\bf v}=(v^1,v^2,v^3)$ is a solution to the Navier--Stokes system (NSE) while $\ep=\pm1$.

Namely, see \cite{KNSS}, the function $u=v_2x_1-v_1x_2\equiv|x'|v_\vartheta$ satisfies the equation
(\ref{NS1}) with $\ep=+1$ (here $v_\vartheta$ is the angular component of the velocity). Next, if
$v_\vartheta=0$, then the function $u=|x'|^{-2}(({\rm rot}({\bf v}))_2x_1-({\rm rot}({\bf v}))_1x_2)$ 
satisfies the equation (\ref{NS1}) with $\ep=-1$.\medskip

Since, by the NSE, $\mbox{div}({\bf v})=0$, it is easy to see that
\begin{equation}\label{div}
\mbox{div}({\bf b})=4\pi\ep\delta_{\Gamma},\qquad\Gamma=\{|x'|=0\}.
\end{equation}
Thus, if $\ep=-1$, then the results of Section 3 are applicable to (\ref{NS1})--(\ref{NS2}). Namely, we are interested in the Liouville theorem.

Note that $\widehat{\bf b}\in L_{q,\infty,loc}(\mathbb R^3\times \mathbb R_-)$ with any $q<2$, and, moreover, satisfies the assumption (\ref{Liouville'}) with $q\in\,]\frac 32,2[$, $\ell=\infty$. Therefore, taking into account
Remark 5, we obtain the following result.\medskip

{\bf Theorem 4.1}. {\it Let ${\bf v}$ be an axisymmetric solution of the Navier--Stokes system in 
$\mathbb R^3\times \mathbb R_-$. Suppose also that ${\bf v}$ satisfies (\ref{Liouville'}) with some $q$ and 
$\ell$ such that $\alpha\equiv\frac 3q+\frac 2\ell-1\in\,[0,1[$. Then any Lipschitz bounded solution of 
(\ref{NS1})--(\ref{NS2}) with $\ep=-1$ in $\mathbb R^n\times\mathbb R_-$ is a constant.}\medskip

{\bf Remark 8}. The assumption (\ref{Liouville'}) is satisfied, for example, if ${\bf v}$ satisfies the estimate
$$|{\bf v}(x',x_3;t)|\le \frac C{|x'|}$$ 
(in this case one can take $q\in\,]\frac 32,2[$, $\ell=\infty$), or
$$|{\bf v}(x',x_3;t)|\le \frac C{(-t)^{\frac 12}}$$ 
(in this case it suffices $q=\infty$, $\ell\in\,]1,2[$).\medskip

To deal with more complicated case $\ep=+1$, we need the following observation.\medskip 

{\bf Remark 9}. The statement of Lemma 3.1 holds true without assumption (\ref{bezdiv}) if $u\le0$ 
in the set ${\cal F}=\mbox{supp}(\mbox{div}({\bf b}))_+$. Similarly, Lemmas 3.2--3.4 and Corollaries 3.1 
(part 2), 3.2, 3.4 hold true if $V\ge k$ in ${\cal F}$. Lemma 3.3 holds true if $V\ge k_0$ in ${\cal F}$.
\medskip

Now we prove the following variant of Corollary 3.3.\medskip

{\bf Lemma 4.2}. {\it Let ${\bf v}$ be an axisymmetric solution of the Navier--Stokes system in $Q_{R}^{2,1}$, and 
let ${\bf v}\in L_{q,\ell}(Q_{R}^{2,1})$ with some $q$ and $\ell$ such that 
$\alpha\equiv\frac 3q+\frac 2\ell-1\in\,[0,1[$.
Let $V$ be a Lipschitz nonnegative supersolution of (\ref{NS1})--(\ref{NS2}) with $\ep=+1$ in $Q_{R}^{2,1}$. If
\begin{equation}\label{Gamma}
V|_{\Gamma\cap Q_{R}^{2,1}}\ge k,\qquad V\le {\mathfrak N}k\quad \mbox{in}\ \ Q_{R}^{2,1} 
\end{equation}
for some $k>0$ and ${\mathfrak N}>1$, then}
\begin{equation}\label{estmin4''}
V\ge \widehat\beta_3 k\qquad\mbox{in}\quad Q_{\frac R2},
\end{equation}
{\it where $\widehat\beta_3$ is completely determined by $q$, $\ell$, ${\mathfrak N}$
and the quantity $\widehat{\cal N}=R^{-\alpha}\|{\bf v}\|_{q,\ell,Q_{R}^{2,1}}$.}

\begin{proof} We put 
$$\widehat{\cal E}_{\varkappa}(t)=\{x\in B_R:\ V(x,t)>\varkappa k\};\qquad
\widehat{\cal E}_{\varkappa}=\{(x;t):\ t\in\, [-R^2,-\frac 34 R^2],\ x\in\widehat{\cal E}_{\varkappa}(t)\}
$$
and claim that 
\begin{equation}\label{measure4}
\mbox{meas}\, (\widehat{\cal E}_{\varkappa})\ge \delta\,\mbox{meas}\,(Q_{R}^{1,\frac 14})
\end{equation}
for some $\varkappa>0$ and $\delta>0$ depending only on $q$, $\ell$, ${\mathfrak N}$ and $\widehat{\cal N}$.

Indeed\footnote{This idea was in a particular case used in \cite{CSTY}.}, by (\ref{div}), we obtain for any Lipschitz test function $\eta\ge0$
\begin{equation}\label{g}
\int\limits_{Q_{R}^{1,\frac 14}(0;-\frac 34 R^2)}\!\!(\partial_tV\eta+ D_iVD_i\eta-b_iVD_i\eta)\,dxdt\ge
4\pi\!\!\int\limits_{\Gamma\cap Q_{R}^{1,\frac 14}(0;-\frac 34 R^2)}\!\!V\eta\,dx_3dt.
\end{equation}
We take $\eta$ such that 
$$ \eta\equiv1\ \ \mbox{in}\ \ Q_{R}^{\frac 12,\frac 18}(0;-\frac{13}{16}R^2);\quad 
\eta\equiv0\ \ \mbox{out of}\ \ Q_{R}^{1,\frac 14}(0;-\frac 34 R^2); 
\qquad |\partial_t\eta|+|D\eta|^2+|\Delta\eta|\le \frac {C}{R^2}.$$ 
Then (\ref{g}) and $V|_{\Gamma\cap Q_{R}^{2,1}}\ge k$ imply
\begin{multline*}
\frac {\pi}2\,kR^3\le \!\!\int\limits_{Q_{R}^{1,\frac 14}(0;-\frac 34 R^2)}\!\!
V\big(|\partial_t\eta|+|\Delta\eta|+|{\bf b}|\cdot|D\eta|\big) dxdt\le\\
\le \frac {C}{R^2}\!\!\int\limits_{Q_{R}^{1,\frac 14}(0;-\frac 34 R^2)}\!\!V\, dxdt+
\frac {C}{R}\!\!\int\limits_{Q_{R}^{1,\frac 14}(0;-\frac 34 R^2)}\!\!V|{\bf v}|\, dxdt+
\frac {C}{R}\!\!\int\limits_{Q_{R}^{1,\frac 14}(0;-\frac 34 R^2)}\!\!\frac V{|x'|}\, dxdt.
\end{multline*}
Splitting the integrals in the right-hand side into integrals over $\widehat{\cal E}_{\varkappa}$ and over 
its complement, we obtain with regard to $V\le {\mathfrak N}k$ and to 
$\frac 1{|x'|}\in L_{\frac 95,\infty}(Q_{R}^{1,\frac 14}(0;-\frac 34 R^2))$
\begin{multline*}
\frac {\pi}2\,kR^3\le \frac{C{\mathfrak N}k}{R^2}\,
\Big[\mbox{meas}\,(\widehat{\cal E}_{\varkappa})+
R^{1+\alpha}\widehat{\cal N}\|{\bf 1}\|_{q',\ell',\,\widehat{\cal E}_{\varkappa}}+
R^{\frac 53}\|{\bf 1}\|_{\frac 94,1,\,\widehat{\cal E}_{\varkappa}}\Big]+\\
+\frac{C\varkappa k}{R^2}\,\Big[\mbox{meas}\,(Q_{R}^{1,\frac 14})+
R^{1+\alpha}\widehat{\cal N}\|{\bf 1}\|_{q',\ell',Q_{R}^{1,\frac 14}}+
R^{\frac 53}\|{\bf 1}\|_{\frac 94,1,Q_{R}^{1,\frac 14}}\Big].
\end{multline*}
The second term in the right-hand side is easily estimated by $C\varkappa kR^3(1+C_{17}(q,\ell)\widehat{\cal N})$. 
Therefore, choosing $\varkappa=\varkappa(q,\ell,\widehat{\cal N})$ sufficiently small, we obtain
\begin{equation}\label{gg}
\frac{C\mathfrak N}{R^5}\Big[\mbox{meas}\,(\widehat{\cal E}_{\varkappa})+
R^{1+\alpha}\widehat{\cal N}\|{\bf 1}\|_{q',\ell',\,\widehat{\cal E}_{\varkappa}}+
R^{\frac 53}\|{\bf 1}\|_{\frac 94,1,\,\widehat{\cal E}_{\varkappa}}\Big]\ge1.
\end{equation}
To estimate the second term in brackets, we rewrite it as follows:
$$\|{\bf 1}\|_{q',\ell',\widehat{\cal E}_{\varkappa}}=\bigg(\int\limits_{-R^2}^{-\frac 34 R^2}\!
\mbox{meas}^{\frac{\ell'}{q'}} (\widehat{\cal E}_{\varkappa}(t))\,dt\bigg)^{\frac 1{\ell'}}.
$$
If $q\le\ell$ then, by the H\"older inequality,
$$\|{\bf 1}\|_{q',\ell',\widehat{\cal E}_{\varkappa}}\le\bigg(\int\limits_{-R^2}^{-\frac 34 R^2}\!
\mbox{meas}\, (\widehat{\cal E}_{\varkappa}(t))\,dt\bigg)^{\frac 1{q'}}\cdot
\Big(\frac 14 R^2\Big)^{\frac 1{\ell'}-\frac 1{q'}}=
\mbox{meas}^{\frac 1{q'}} (\widehat{\cal E}_{\varkappa})\cdot\Big(\frac 14 R^2\Big)^{\frac 1q-\frac 1{\ell}}.
$$
In the opposite case, since $\mbox{meas}\, (\widehat{\cal E}_{\varkappa}(t))\le \mbox{meas}\, (B_R)$, we obtain
\begin{multline*}
\|{\bf 1}\|_{q',\ell',\widehat{\cal E}_{\varkappa}}=\bigg(\int\limits_{-R^2}^{-\frac 34 R^2}\!
\bigg(\frac{\mbox{meas}\, (\widehat{\cal  E}_{\varkappa}(t))}{\mbox{meas}\,(B_R)}\bigg)
^{\frac{\ell'}{q'}}dt\bigg)^{\frac 1{\ell'}}\cdot \mbox{meas}^{\frac 1{q'}}(B_R)\le\\
\le \bigg(\int\limits_{-R^2}^{-\frac 34 R^2}\!
\frac{\mbox{meas}\, (\widehat{\cal  E}_{\varkappa}(t))}{\mbox{meas}\,(B_R)}
\,dt\bigg)^{\frac 1{\ell'}}\cdot \mbox{meas}^{\frac 1{q'}}(B_R)=
\mbox{meas}^{\frac 1{\ell'}} (\widehat{\cal E}_{\varkappa})\cdot(4\pi R^3)^{\frac 1{\ell}-\frac 1q}.
\end{multline*}
Similarly we estimate the third term in brackets in (\ref{gg}). Thus, we obtain the inequality for 
${\cal A}=\frac 1{R^5}\,\mbox{meas}\, (\widehat{\cal  E}_{\varkappa})$:
$${\cal A}+{\cal A}^{\frac 49}+\widehat{\cal N}{\cal A}^{1-\max\{\frac 1q,\frac 1{\ell}\}}\ge \frac 1{C{\mathfrak N}},
$$
and (\ref{measure4}) follows.\medskip

The inequality (\ref{measure4}) provides 
$$\mbox{meas}\, (\widehat{\cal E}_{\varkappa}(\bar t))\ge \delta\,\mbox{meas}\,(B_{R})
$$
for some $\bar t\in\,[-R^2, -\frac 34 R^2]$, and Corollary 3.3 ensures (\ref{estmin4''}) with
$\widehat\beta_3=\varkappa\cdot\beta_3(\delta, 3, 1, q, \ell, \widehat{\cal N})$.
\end{proof}

{\bf Theorem 4.3}. {\it Let ${\bf v}$ be an axisymmetric solution of the Navier--Stokes system in 
$\mathbb R^3\times \mathbb R_-$. Suppose also that ${\bf v}$ satisfies (\ref{Liouville'}) with some $q$ and 
$\ell$ such that $\alpha\equiv\frac 3q+\frac 2\ell-1\in\,[0,1[$. Let $u$ be a Lipschitz bounded solution of 
(\ref{NS1})--(\ref{NS2}) with $\ep=+1$ in $\mathbb R^n\times\mathbb R_-$. If $u|_{\Gamma}=const$, then 
$u\equiv const$.}\medskip

\begin{proof} Given $R$, we set $k=\frac 12\underset{Q_{2R}}{\mbox{osc}}\ u$ and consider two functions $V_1=u-\inf\limits_{Q_{2R}}u$ and $V_2=\sup\limits_{Q_{2R}}u-u$. At least one of them satisfies 
(\ref{Gamma}) with ${\mathfrak N}=2$. Therefore, Lemma 4.2 gives for this function the estimate
(\ref{estmin4''}), which implies (\ref{osc1}) with 
$\varkappa_1=1-\frac 12\widehat\beta_3(q,\ell,2,\widehat{\cal N})$. Iteration of this 
inequality with respect to a suitable sequence $R_m\to\infty$ completes the proof.
\end{proof}


\begin{thebibliography}{AFT}

\bibitem[DG]{DG}
E. De Giorgi, {\em Sulla differenziabilit\`a e l'analiticit\`a delle estremali degli integrali multipli
regolari}, Mem. Accad. Sci. Torino. Cl. Sci. Fis. Mat. Nat. {\bf 3} (1957), 25--43. 

\bibitem[M]{M}
C.B. Morrey, Jr., {\em Second order elliptic equations in several variables and H\"older continuity},
Math. Z. {\bf 72} (1959/1960), 146--164.

\bibitem[Na]{Na}
J. Nash, {em Parabolic equations}, Proc. Nat. Acad. Sci. USA, {\bf 43} (1957), 754--758. 

\bibitem[LU1]{LU1}
O.A. Lady\v{z}enskaja, N.N. Ural'ceva, {\em A boundary-value problem for linear and quasi-linear 
parabolic equations}, Dokl. Akad. Nauk SSSR 139 1961 544--547 (Russian).

\bibitem[Mo1]{Mo1}
J. Moser, {\em On Harnack's theorem for elliptic differential equations}, Comm. Pure Appl. Math. {\bf 14} 
(1961), 577--591. 

\bibitem[Mo2]{Mo2}
J. Moser, {\em A Harnack inequality for parabolic differential equations}, Comm. Pure Appl. Math. {\bf 17}
(1964), 101--134. Corrected in: Comm. Pure Appl. Math. {\bf 20} (1967), 231--236. 

\bibitem[NU]{NU}
A.I. Nazarov, N.N. Ural'tseva, {\em Qualitative properties of solutions 
to elliptic and parabolic equations with unbounded lower-order coefficients},
SPbMS El. Prepr. Archive. N~2009-05. 6p.

\bibitem[Z]{Z}
Qi S. Zhang, {\em A strong regularity result for parabolic equations},
Comm. Math. Phys. {\bf 244} (2004), N2, 245--260. 

\bibitem[KNSS]{KNSS}
G. Koch, N. Nadirashvili, G. Seregin, V. Sverak, {\em Liouville theorems for
the Navier--Stokes equations and applications}, Acta Math. {\bf 203} (2009), N1,
83--105.

\bibitem[CSTY]{CSTY}
C.-C. Chen, R.M. Strain, T.-P. Tsai, H.-T. Yau,{\em Lower bound on the blow-up rate 
of the axisymmetric Navier--Stokes equations}, I, 
Int. Math. Res. Not., {\bf 8} (2008), Art. ID rnn016, 31 pp.; II,
Comm. PDE, {\bf 34} (2009), N1-3, 203--232.

\bibitem[SSSZ]{SSSZ}
G. Seregin, L. Silvestre, V. Sverak, A. Zlatos, {\em On divergence-free drifts},
preprint arXiv:1010.6025v1

\bibitem[Tru]{Tru}
N.S. Trudinger, {\em Linear elliptic operators with measurable coefficients}, Ann. Scuola Norm. Sup. Pisa,
{\bf 27} (1973), 265--308. 

\bibitem[Li1]{Li1}
G.M. Lieberman, {\em Second order parabolic differential equations}, World Sci. Publishing Co., 
NJ, 1996. 

\bibitem[LU2]{LU2}
O.A. Ladyzhenskaya, N.N. Uraltseva, {\em Linear and quasilinear elliptic equations},
2nd ed., Moscow, Nauka, 1973  (Russian); English transl. of the 1st ed.: 
Academic Press, NY, 1968.

\bibitem[LSU]{LSU}
Ladyzhenskaja, O.A.; Solonnikov, V.A.; Ural'tseva, N.N., {\em
Linear and quasi-linear equations of parabolic type}, Moscow,
Nauka, 1968 (Russian); English transl.: Translations of
Mathematical Monographs, {\bf 23}, AMS, Providence, R.I., 1967.

\bibitem[LL]{LL} 
E. Lieb, M. Loss, {\em Analysis}, Grad. Studies in Math. {\bf 14}, 2nd ed., AMS, 2001.

\bibitem[S2]{S2}
M.V. Safonov, {\em Mean value theorems and Harnack inequalities for second-order parabolic equations}, 
Nonlinear problems in mathematical physics and related topics, II, Int. Math. Ser., {\bf 2}, 
Kluwer/Plenum, NY, 2002, 329--352. 

\bibitem[MV]{MV}
V.G. Maz'ya, I.E. Verbitsky,  {\em Form boundedness of the general second order differential operator},
Comm. Pure Appl. Math. {\bf 59} (2006), N9, 1286--1329. 

\bibitem[Tro]{Tro}
G.M. Troianiello, {\em Elliptic Differential Equations and Obstacle Problems},  NY, 1987. 

\bibitem[BIN]{BIN}
Besov, O.V.; Il'in, V.P.; Nikol'skii, S.M., {\em Integral
representations of functions and imbedding theorems}, ed.2,
Moscow, Nauka, 1996 (Russian); English transl. of the 1st ed.:
V. 1.  Scripta Series in Mathematics. Edited by Mitchell H.
Taibleson. V. H. Winston \& Sons, Washington, D.C.; Halsted Press
[John Wiley \& Sons], New York-Toronto, Ont.-London, 1978.

\bibitem[S3]{S3}
M.V. Safonov, {\em Non-divergence elliptic equations of second order with unbounded drift},
AMS Transl. Series 2. {\bf 229} (2010), 211--232.

\end{thebibliography}
\end{document}